\documentclass[reqno]{gokova}
\usepackage{epsfig,graphicx}
\usepackage{amssymb}

\newcommand{\ex}{\mathfrak}
\newcommand{\TC}{{}^{\text{log}}_{\phantom{t}\mathbb T}C}
\newcommand{\logT}{{}^\text{log} T}
\newcommand{\logC}{{}^\text{log}C}
\providecommand{\fp}[2]{{}_{\hspace{3pt}#1\hspace{-2pt}}\times_{#2}}

\newcommand{\dbar}{\bar{\partial}}  
\newcommand{\exploded}{exploded $\mathbb T$ }
\providecommand{\abs}[1]{\left\lvert #1\right\rvert}

\begin{document}
\def\E{\ifmmode{\mathbb E}\else{$\mathbb E$}\fi} 
\def\N{\ifmmode{\mathbb N}\else{$\mathbb N$}\fi} 
\def\R{\ifmmode{\mathbb R}\else{$\mathbb R$}\fi} 
\def\Q{\ifmmode{\mathbb Q}\else{$\mathbb Q$}\fi} 
\def\C{\ifmmode{\mathbb C}\else{$\mathbb C$}\fi} 
\def\H{\ifmmode{\mathbb H}\else{$\mathbb H$}\fi} 
\def\Z{\ifmmode{\mathbb Z}\else{$\mathbb Z$}\fi} 
\def\P{\ifmmode{\mathbb P}\else{$\mathbb P$}\fi} 
\def\T{\ifmmode{\mathbb T}\else{$\mathbb T$}\fi} 
\def\SS{\ifmmode{\mathbb S}\else{$\mathbb S$}\fi} 
\def\DD{\ifmmode{\mathbb D}\else{$\mathbb D$}\fi} 

\renewcommand{\a}{\alpha}
\renewcommand{\b}{\beta}
\renewcommand{\d}{\delta}
\newcommand{\D}{\Delta}
\newcommand{\e}{\varepsilon}
\newcommand{\g}{\gamma}
\newcommand{\G}{\Gamma}
\newcommand{\la}{\lambda}
\newcommand{\La}{\Lambda}
\newcommand{\n}{\nabla}
\newcommand{\var}{\varphi}
\newcommand{\s}{\sigma}
\newcommand{\Sig}{\Sigma}
\renewcommand{\t}{\tau}
\renewcommand{\th}{\theta}
\renewcommand{\O}{\Omega}
\renewcommand{\o}{\omega}
\newcommand{\z}{\zeta}

\newcommand{\ben}{\begin{enumerate}}
\newcommand{\een}{\end{enumerate}}
\newcommand{\be}{\begin{equation}}
\newcommand{\ee}{\end{equation}}
\newcommand{\bea}{\begin{eqnarray}}
\newcommand{\eea}{\end{eqnarray}}
\newcommand{\bc}{\begin{center}}
\newcommand{\ec}{\end{center}}

\newtheorem{thm}{Theorem}[section]
\newtheorem{cor}[thm]{Corollary}
\newtheorem{lem}[thm]{Lemma}
\newtheorem{prop}[thm]{Proposition}
\newtheorem{ax}{Axiom}
\newtheorem{conj}[thm]{Conjecture}

\theoremstyle{definition}
\newtheorem{defn}{Definition}[section]

\theoremstyle{remark}
\newtheorem{rem}{\rm\bfseries{Remark}}[section]
\newtheorem*{notation}{Notation}

\newtheorem{ques}{\rm\bfseries{Question}}[section]
\newtheorem{cons}[rem]{\rm\bfseries{Construction}}
\newtheorem{exm}[rem]{\rm\bfseries{Example}}



\setcounter{page}{1}
\volume{13}

\title[Exploded fibrations]{Exploded fibrations}
\author[PARKER]{Brett Parker}

\thanks{This research was partially supported by money from the NSF grant 0244663.}

\address{Department of Mathematics, MIT, Cambridge, MA 02139-4307, USA}
\email{parker@math.mit.edu}

\begin{abstract}

Initiated by Gromov in \cite{gromov}, the study of holomorphic curves in symplectic manifolds has been a powerfull tool in symplectic topology, however the moduli space of holomorphic curves is often very difficult to find. A common technique is to study the limiting behavior of holomorphic curves in a degenerating family of complex structures which corresponds to a kind of adiabatic limit. The category of exploded fibrations is an extension of the smooth category in which some of these degenerations can be described as smooth families.

The first part of this paper is devoted to defining exploded fibrations and a slightly more specialized category of exploded $\mathbb T$ fibrations. In section \ref{examples} are some examples of holomorphic curves in exploded $\mathbb T$ fibrations, including a brief discussion of the relationship between tropical geometry and exploded $\mathbb T$ fibrations. In section \ref{fiber product}, it is shown that exploded fibrations have a good intersection theory. In section \ref{perturbation theory}, the perturbation theory of holomorphic curves in \exploded fibrations is sketched.

\end{abstract}

\maketitle

\section{Introduction}

A symplectic manifold is a manifold $M^{2n}$ with a closed, maximally nondegenerate two form $\omega$, called the symplectic form. Every manifold of this type has local coordinates $(x,y)\in \mathbb R^n\times \mathbb R^n$ in which the symplectic form looks like $\sum dx_i\wedge dy_i$. For this reason symplectic manifolds have no local invariants, and the study of symplectic manifolds is called symplectic topology.

One of the most powerfull tools in symplectic topology is the study of holomorphic curves. Given a symplectic manifold $(M,\omega)$, there is a contractible choice of almost complex structure $J$ on $M$ which is tamed by $\omega$ in the sense that $\omega(v,Jv)>0$ for any nonzero vector $v$. With such a choice of $J$, a holomorphic curve is a map 

\[f:(S,j)\longrightarrow (M,J)\] 
from a Riemann surface $S$ with a complex structure $j$ so that $df\circ j=J\circ df$. (These are sometimes called pseudoholomorphic curves because $(M,J)$ is not a complex manifold.) The energy of a holomorphic curve is $\int_S f^*\omega$. In \cite{gromov}, Gromov proved that the moduli space of  holomorphic curves with bounded genus and energy in a compact symplectic manifold  is compact. The moduli space should also be `smooth' in some sense, and an isotopy of $J$ gives a compact cobordism of moduli spaces.  This allows the definition of Gromov Witten invariants of a symplectic manifold which involve counting the number of holomorphic curves with some constraints.

\

The sense in which the moduli space of holomorphic curves is `smooth' is not the usual definition. Usually a smooth family of maps corresponds to a pair of smooth maps
\[F\longleftarrow C\longrightarrow M\] 
where $C\longrightarrow F$ is a surjective submersion. The individual maps in this family are the maps of the fibers of $C\longrightarrow F$ into $M$. 

The moduli space of holomorphic curves includes families of holomorphic curves in which a bubble forms

\
\begin{figure}[htb]
\includegraphics[width=4 in]{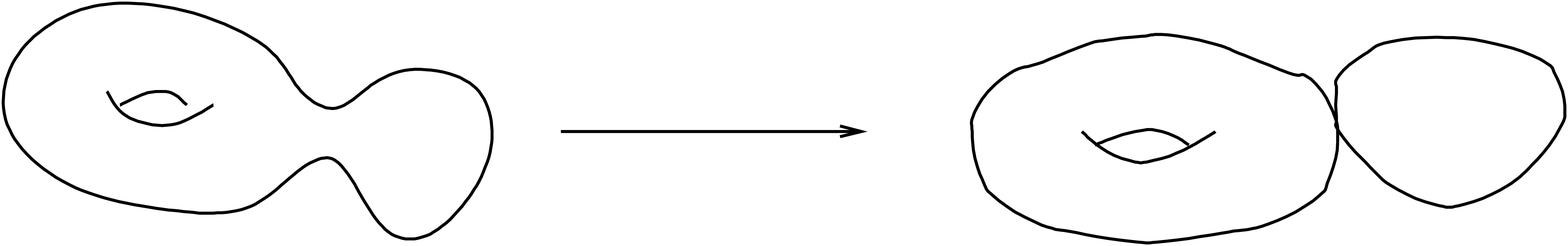}
\end{figure}

This change in the topology of the domain can not take place in a connected smooth family of maps, so if the behavior above is to be considered `smooth', we need a new definition. This is one indication that the category of smooth manifolds might not be the natural category to describe the theory of  holomorphic curves.

\

A second reason to try to find an  extension of the smooth category with a good theory of holomorphic curves is that holomorphic curves are in general difficult to find in the non algebraic setting. Many techniques for finding holomorphic curve invariants involve degenerating the almost complex structure $J$ in some way so that holomorphic curves become easier to find in the limit. We will now see two types of degenerations which can be considered as cutting a symplectic manifold into smaller, simpler pieces in order to compute holomorphic curve invariants. These degenerations do not have a smooth limit in the category of manifolds, however they can be considered as `smooth' families of exploded fibrations.

\

One example similar to the formation of the bubble shown above is the degeneration of almost complex structure used in symplectic field theory to break a symplectic manifold apart along a hypersurface. This is described in \cite{sft}. 

Suppose we have a hypersurface $S\subset M$ with a collar neighborhood equal to $(0,1)\times S$. Denote by $t$ the coordinate for $(0,1)$, by $X_t$ the Hamiltonian vectorfield generated by $t$ so that  $\omega(X_t,\cdot)=dt$. Suppose also that the vectorfield $X_t$ on $S$ and $\omega(\frac\partial {\partial t},\cdot)$  are both independent of $t$. We can choose an almost complex structure $J$ tamed by $\omega$ which is independent of $t$, so that $JX_t=\frac\partial {\partial t}$, and which preserves the kernel of $\omega(\frac\partial{\partial t},\cdot )$ restricted to $S$.

\

\begin{figure}[htb]

\begin{picture}(0,0)%
\includegraphics{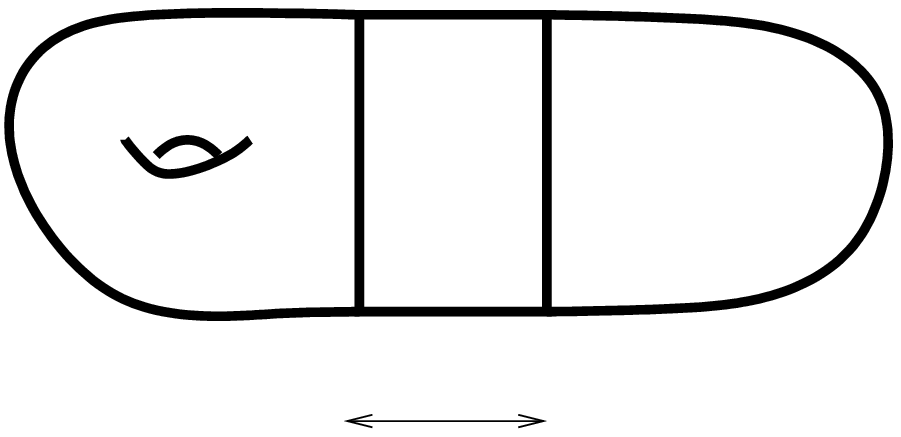}%
\end{picture}%
\setlength{\unitlength}{3947sp}%
\begingroup\makeatletter\ifx\SetFigFont\undefined%
\gdef\SetFigFont#1#2#3#4#5{%
  \reset@font\fontsize{#1}{#2pt}%
  \fontfamily{#3}\fontseries{#4}\fontshape{#5}%
  \selectfont}%
\fi\endgroup%
\begin{picture}(4307,2277)(1951,-5735)
\put(3901,-5686){\makebox(0,0)[lb]{\smash{{\SetFigFont{12}{14.4}{\rmdefault}{\mddefault}{\updefault}{$t$}%
}}}}
\put(3708,-4253){\makebox(0,0)[lb]{\smash{{\SetFigFont{12}{14.4}{\rmdefault}{\mddefault}{\updefault}$S\times(0,1)$}%
}}}
\end{picture}%

\end{figure} 

 We can then choose a degenerating family of almost complex structures $J^\epsilon$ on $M$ as follows. Keep $J^\epsilon$ constant on the kernel of $\omega(\frac\partial{\partial t},\cdot)$ restricted to $S$, and have 
\[J^\epsilon X_t={\frac {dt}{dT^\epsilon}}\frac \partial{\partial t}\]
where $T^\epsilon(t)$ is some family of smooth monotone increasing functions on $[0,1]$  with $\frac {dT^\epsilon}{dt}=1$ in a neighborhood of $0$ and $1$ and so that $T^\epsilon(0)=\log \epsilon$ and $T(1)=1$. From the perspective of the almost complex structure, $J^\epsilon$ has the effect of replacing our neighborhood $(0,1)\times S$ with $(\log\epsilon,1)\times S$, where the almost complex structure on this lengthened cylinder with coordinate $T$ is just the symmetric extension of the old $J$.

\begin{figure}[htb]

\begin{picture}(0,0)%
\includegraphics{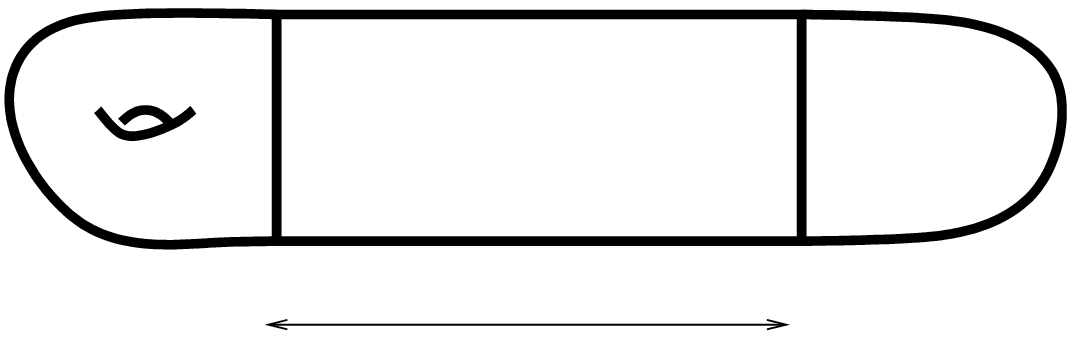}%
\end{picture}%
\setlength{\unitlength}{3947sp}%
\begingroup\makeatletter\ifx\SetFigFont\undefined%
\gdef\SetFigFont#1#2#3#4#5{%
  \reset@font\fontsize{#1}{#2pt}%
  \fontfamily{#3}\fontseries{#4}\fontshape{#5}%
  \selectfont}%
\fi\endgroup%
\begin{picture}(5141,1757)(1155,-3424)
\put(3570,-3379){\makebox(0,0)[lb]{\smash{{\SetFigFont{9}{10.8}{\rmdefault}{\mddefault}{\updefault}$T^\epsilon$}%
}}}
\put(3018,-2183){\makebox(0,0)[lb]{\smash{{\SetFigFont{12}{14.4}{\rmdefault}{\mddefault}{\updefault}$S\times (\log\epsilon,1)$}%
}}}
\end{picture}%
\end{figure}

A short calculation shows that $J^\epsilon$ is still tamed by $\omega$. Note that we have a great amount of flexibility in how the almost complex manifold $(M,J^\epsilon)$ wears the symplectic form $\omega$. Different choices of $T^\epsilon$ will concentrate $\omega$ in different regions of the cylinder. This is important for being able to tame holomorphic curves in the limit as $\epsilon\rightarrow 0$.

Suppose that this hypersurface $S$ separates $M$ into $M_0$ and $M_1$. In the limit as $\epsilon\rightarrow 0$, our almost complex manifold becomes modelled on the manifold with cylindrical ends $\tilde M_0$ obtained by gluing   $(0,\infty)\times S$ to $M_0$, the manifold $\tilde M_1$ obtained by gluing  $(-\infty,1)\times S$ to $M_1$, or the cylinder $\mathbb R\times S$. In some sense, the `limit' of this family of almost complex manifolds should include an infinite number of these cylindrical pieces, because one can construct a sequence of maps of $\mathbb Z$ to $(M,J^\epsilon)$ so that in the limit, the model around the image of each point is a cylinder which contains none of the other points.

\begin{figure}[htb]
 \begin{picture}(0,0)%
\includegraphics{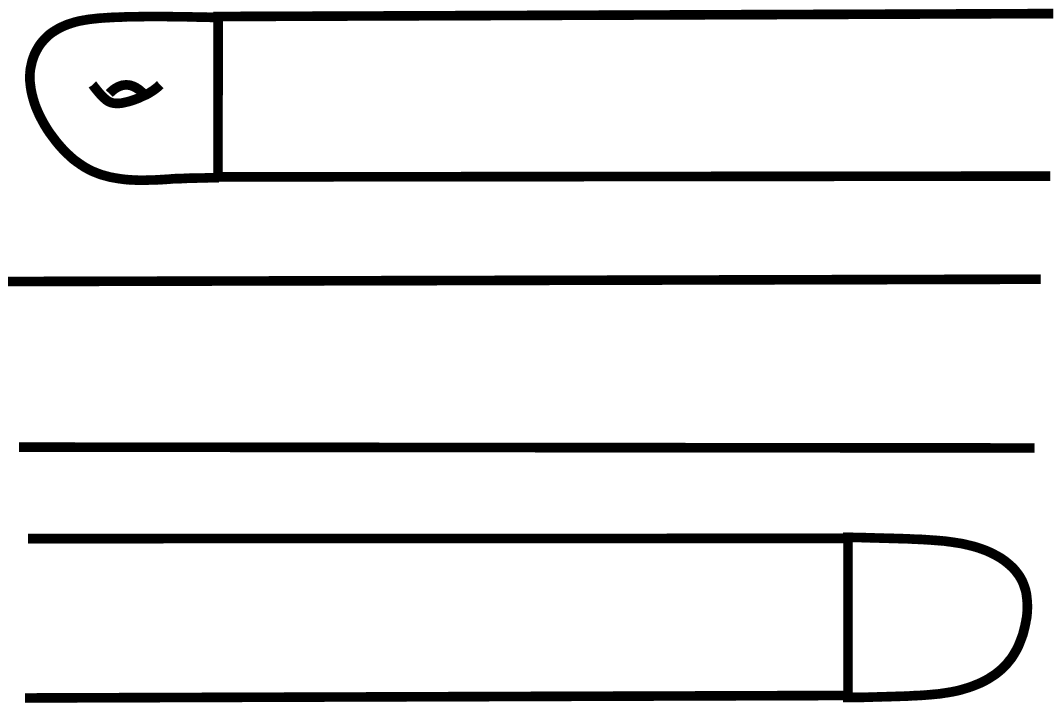}%
\end{picture}%
\setlength{\unitlength}{3947sp}%
\begingroup\makeatletter\ifx\SetFigFont\undefined%
\gdef\SetFigFont#1#2#3#4#5{%
  \reset@font\fontsize{#1}{#2pt}%
  \fontfamily{#3}\fontseries{#4}\fontshape{#5}%
  \selectfont}%
\fi\endgroup%
\begin{picture}(5104,3373)(2306,-4665)
\put(3998,-2994){\makebox(0,0)[lb]{\smash{{\SetFigFont{9}{10.8}{\rmdefault}{\mddefault}{\updefault}{$\mathbb R\times S$}%
}}}}
\put(3941,-1774){\makebox(0,0)[lb]{\smash{{\SetFigFont{9}{10.8}{\rmdefault}{\mddefault}{\updefault}{$\tilde M_0$}%
}}}}
\put(4136,-4274){\makebox(0,0)[lb]{\smash{{\SetFigFont{9}{10.8}{\rmdefault}{\mddefault}{\updefault}{$\tilde M_1$}%
}}}}
\end{picture}%

\end{figure}

We will be able to view this family $(M,J^\epsilon)$ as part of a connected `smooth' family of exploded fibrations which contains an object that has a base equal to  the interval $[0,1]$ and fibers equal to $\tilde M_0$ over $0$, $\tilde M_1$ over $1$ and a copy of $\mathbb R\times S$ over each point in the interior.

\begin{figure}[htb]
 \begin{picture}(0,0)%
\includegraphics{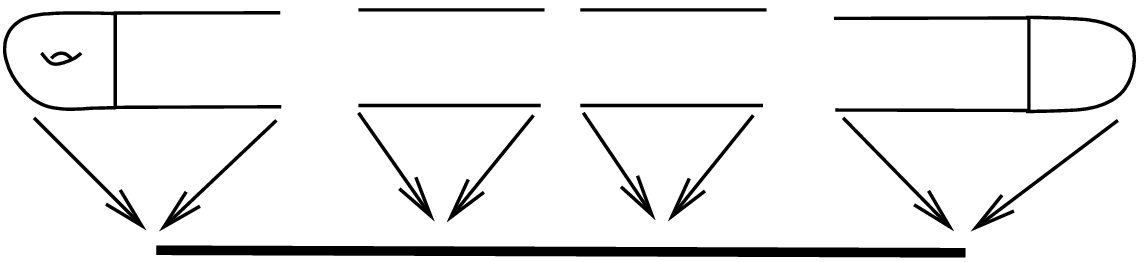}%
\end{picture}%
\setlength{\unitlength}{3947sp}%
\begingroup\makeatletter\ifx\SetFigFont\undefined%
\gdef\SetFigFont#1#2#3#4#5{%
  \reset@font\fontsize{#1}{#2pt}%
  \fontfamily{#3}\fontseries{#4}\fontshape{#5}%
  \selectfont}%
\fi\endgroup%
\begin{picture}(5467,1232)(858,-5884)
\put(2804,-4891){\makebox(0,0)[lb]{\smash{{\SetFigFont{6}{7.2}{\rmdefault}{\mddefault}{\updefault}{$\mathbb R\times S$}%
}}}}
\put(3870,-4891){\makebox(0,0)[lb]{\smash{{\SetFigFont{6}{7.2}{\rmdefault}{\mddefault}{\updefault}{$\mathbb R\times S$}%
}}}}
\put(1754,-4939){\makebox(0,0)[lb]{\smash{{\SetFigFont{6}{7.2}{\rmdefault}{\mddefault}{\updefault}{$\tilde M_0$}%
}}}}
\put(4979,-4959){\makebox(0,0)[lb]{\smash{{\SetFigFont{6}{7.2}{\rmdefault}{\mddefault}{\updefault}{$\tilde M_1$}%
}}}}
\end{picture}%

\end{figure}

In section \ref{log}, we will describe the exact structure on these fibers, in section \ref{base}, we will describe the structure on the base, and in section \ref{exfib}, we will say exactly how these fit together. 
The category of exploded  fibrations has a well defined product, so we will also be able to deal with degenerations which look locally like products of degenerations of this type.

If the dynamics of $X_t$ satisfy a certain condition (which holds generically), the limit of $J^\epsilon$ holomorphic curves under this degeneration look like a number of holomorphic curves mapping into the different model fibers which are asymptotic to cylinders over orbits of $X_t$. The holomorphic curves inside all but a finite number of the cylindrical models over the interior of the interval consist just of cylinders over orbits of $X_t$ which coincide with the asymptotics of the next nontrivial holomorphic curves in models to the left or right. 

\

\begin{figure}[htb]
\includegraphics{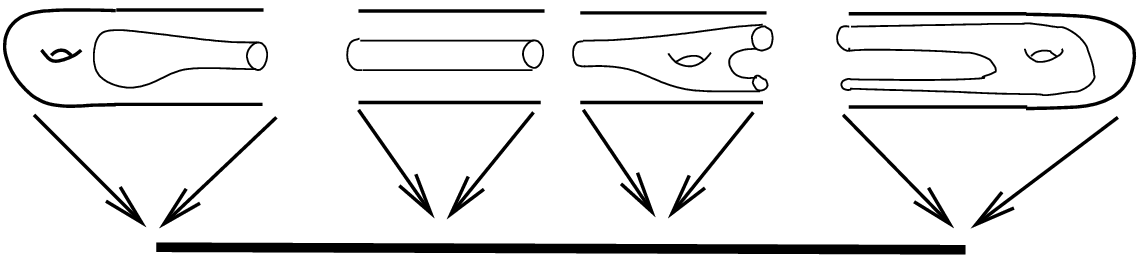}
\end{figure}

\

Such holomorphic curves in exploded fibrations are morphisms in the exploded category with a slightly weaker `smooth' structure. A special type of exploded fibration where holomorphic curves are better behaved is called an exploded $\mathbb T$ fibration. The degeneration above would give an exploded $\mathbb T$ fibration if the action of $X_t$ was a free circle action. We will define the category of exploded $\mathbb T$ fibrations in section \ref{etf}, and sketch the theory of holomorphic curves in exploded $\mathbb T$ fibrations in section \ref{examples}.

\
 
For a second example, consider $\mathbb CP^2$ with the standard coordinates $[z_0,z_1,z_2]$, and an action of $\mathbb T^2$ by $[z_0,z_1,z_2]\mapsto[z_0,e^{\theta_1i}z_1,e^{\theta_2i}z_2]$. We can  choose a symplectic form taming the complex structure on $\mathbb CP^2$ which is preserved by this $\mathbb T^2$ action. Such an $\omega$ is of the form $dh_1\wedge d \theta_1+dh_2\wedge d\theta_2$ where $h_1$ and $h_2$ are $\mathbb T^2$ invariant functions. The functions $h_1$ and $h_2$ are the Hamiltonian functions generating the action of $\theta_1$ and $\theta_2$. The image of $(h_1,h_2)$ is called the moment polytope. In the case of $\mathbb CP^2$, we can choose this to be the triangle  $\{0\leq h_1,0\leq h_2, 5\geq h_1+h_2\}$. There is a large amount of flexibility in how we choose $h_1$ and $h_2$. In particular, we can make the following choices 

\begin{enumerate}\item $h_1=\log\abs{\frac {z_1}{z_0}}$ and $h_2=h_2\left(\abs{\frac {z_2}{z_0}}\right)$ on $\{1\leq h_1\leq 3, h_1+h_2\leq 4\}$ 

\item $h_2=\log\abs{\frac{z_2}{z_0}}$ and $h_1=h_1\left(\abs{\frac{z_1}{z_0}}\right)$ on $\{1\leq h_2\leq 3, h_1+h_2\leq 4\}$
\item $h_1-h_2=\log\abs{\frac {z_1}{z_2}}$ and $h_1+h_2=(h_1+h_2)\left(\abs{\frac{z_0^2}{z_1z_2}}\right)$ \\ on $\{-2\leq h_1-h_2\leq 2, h_1+h_2\geq 4\}$ 
\end{enumerate}

This gives a complex structure with the symmetry over the moment polytope shown in the figure below.

\

\begin{figure}[htb]
 \begin{picture}(0,0)%
\includegraphics{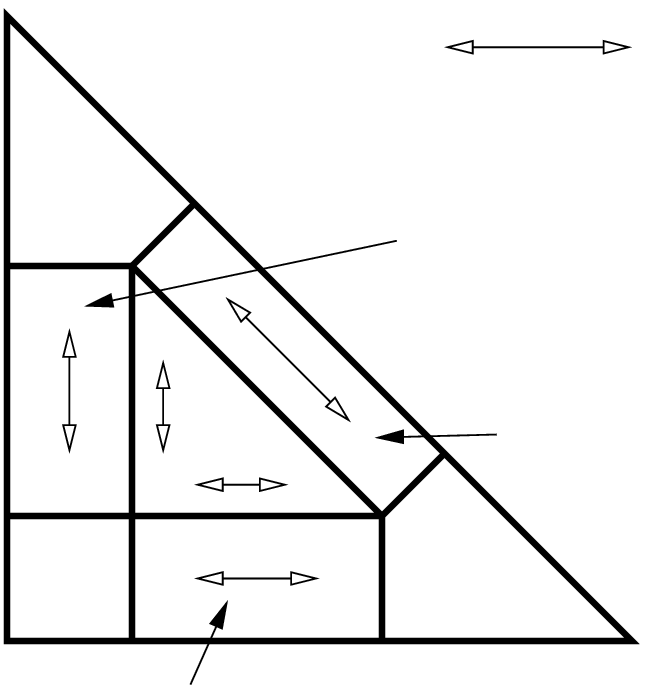}%
\end{picture}%
\setlength{\unitlength}{3947sp}%
\begingroup\makeatletter\ifx\SetFigFont\undefined%
\gdef\SetFigFont#1#2#3#4#5{%
  \reset@font\fontsize{#1}{#2pt}%
  \fontfamily{#3}\fontseries{#4}\fontshape{#5}%
  \selectfont}%
\fi\endgroup%
\begin{picture}(3267,3556)(1168,-5084)
\put(1351,-4261){\makebox(0,0)[lb]{\smash{{\SetFigFont{12}{14.4}{\rmdefault}{\mddefault}{\updefault}{$M_1$}%
}}}}
\put(1351,-2311){\makebox(0,0)[lb]{\smash{{\SetFigFont{12}{14.4}{\rmdefault}{\mddefault}{\updefault}{$M_2$}%
}}}}
\put(3301,-4261){\makebox(0,0)[lb]{\smash{{\SetFigFont{12}{14.4}{\rmdefault}{\mddefault}{\updefault}{$M_3$}%
}}}}
\put(2251,-2011){\makebox(0,0)[lb]{\smash{{\SetFigFont{12}{14.4}{\rmdefault}{\mddefault}{\updefault}{ translation symmetry in $J$}%
}}}}
\put(3376,-3436){\makebox(0,0)[lb]{\smash{{\SetFigFont{12}{14.4}{\rmdefault}{\mddefault}{\updefault}{$S_{2,3}\times \{1\leq\abs z\leq e^4\}$}%
}}}}
\put(2691,-2465){\makebox(0,0)[lb]{\smash{{\SetFigFont{12}{14.4}{\rmdefault}{\mddefault}{\updefault}{$S_{1,2}\times\{1\leq\abs z\leq e^3\}$}%
}}}}
\put(1621,-5026){\makebox(0,0)[lb]{\smash{{\SetFigFont{12}{14.4}{\rmdefault}{\mddefault}{\updefault}{$S_{1,3}\times\{1\leq\abs z\leq e^3\}$}%
}}}}
\end{picture}%

\end{figure}

\

\

We will now draw a degeneration of complex structure that could also be applied if we glued  anything into the three regions $M_1$ $M_2$ and $M_3$ in the above picture. As in the previous example, there is a large amount of flexibility in how this new  complex manifold wears the symplectic structure. We choose some way which doesn't change in the regions $M_1$ $M_2$ and $M_3$.

\newpage

\begin{figure}[htb]
 \begin{picture}(0,0)%
\includegraphics{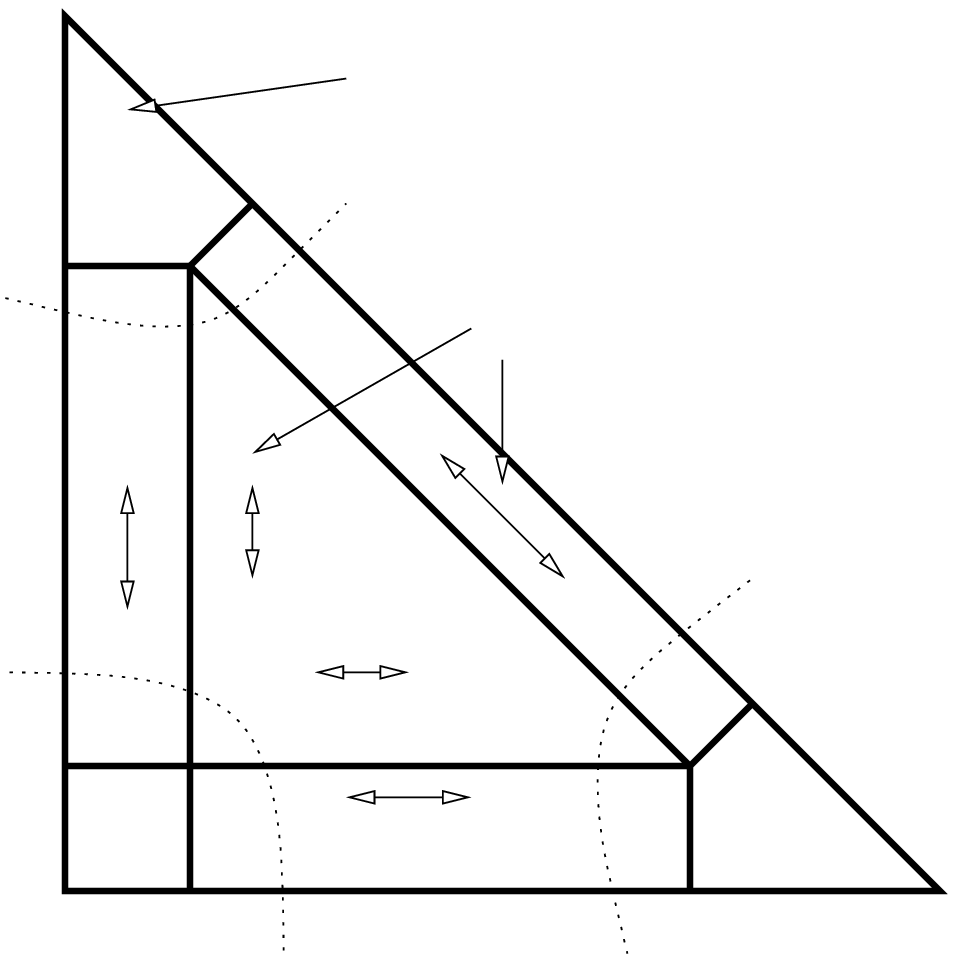}%
\end{picture}%
\setlength{\unitlength}{3947sp}%
\begingroup\makeatletter\ifx\SetFigFont\undefined%
\gdef\SetFigFont#1#2#3#4#5{%
  \reset@font\fontsize{#1}{#2pt}%
  \fontfamily{#3}\fontseries{#4}\fontshape{#5}%
  \selectfont}%
\fi\endgroup%
\begin{picture}(4934,4545)(889,-6073)
\put(1351,-2311){\makebox(0,0)[lb]{\smash{{\SetFigFont{12}{14.4}{\rmdefault}{\mddefault}{\updefault}{ $M_2$}%
}}}}
\put(1351,-5461){\makebox(0,0)[lb]{\smash{{\SetFigFont{12}{14.4}{\rmdefault}{\mddefault}{\updefault}{ $M_1$}%
}}}}
\put(4351,-5461){\makebox(0,0)[lb]{\smash{{\SetFigFont{12}{14.4}{\rmdefault}{\mddefault}{\updefault}{ $M_3$}%
}}}}
\put(2701,-1861){\makebox(0,0)[lb]{\smash{{\SetFigFont{12}{14.4}{\rmdefault}{\mddefault}{\updefault}{ Use old $J$ here}%
}}}}
\put(3301,-3061){\makebox(0,0)[lb]{\smash{{\SetFigFont{12}{14.4}{\rmdefault}{\mddefault}{\updefault}{ Complete $J$ symmetrically here}%
}}}}
\put(2401,-5611){\makebox(0,0)[lb]{\smash{{\SetFigFont{12}{14.4}{\rmdefault}{\mddefault}{\updefault}{ $S_{1,3}\times\{1\leq\abs z\leq \frac{e^3}{\epsilon}\}$}%
}}}}
\end{picture}%

\end{figure}

Imagine rescaling this triangle to keep it a constant size as we change $J$. In the limit, this will give a base for a limiting exploded $\mathbb T$ fibration. The fiber over each point in this base will be the local model for what the complex structure looks like around this point in the limit.

\newpage

\begin{figure}[htb]
 \begin{picture}(0,0)%
\includegraphics{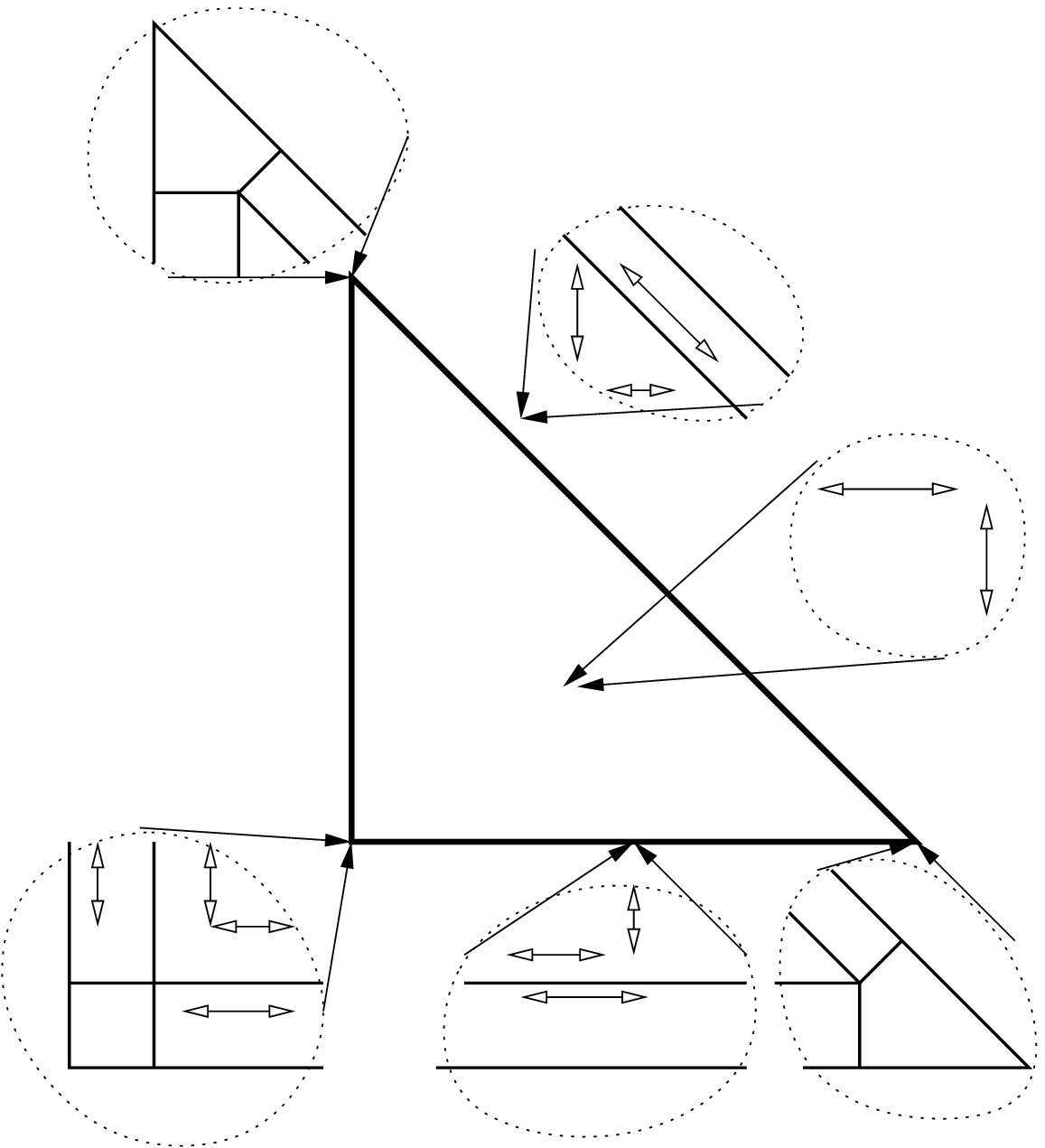}%
\end{picture}%
\setlength{\unitlength}{3947sp}%
\begingroup\makeatletter\ifx\SetFigFont\undefined%
\gdef\SetFigFont#1#2#3#4#5{%
  \reset@font\fontsize{#1}{#2pt}%
  \fontfamily{#3}\fontseries{#4}\fontshape{#5}%
  \selectfont}%
\fi\endgroup%
\begin{picture}(5520,6070)(-668,-6188)
\put(-299,-5611){\makebox(0,0)[lb]{\smash{{\SetFigFont{12}{14.4}{\rmdefault}{\mddefault}{\updefault}{ $M_1$}%
}}}}
\put(226,-811){\makebox(0,0)[lb]{\smash{{\SetFigFont{12}{14.4}{\rmdefault}{\mddefault}{\updefault}{ $M_2$}%
}}}}
\put(3976,-5611){\makebox(0,0)[lb]{\smash{{\SetFigFont{12}{14.4}{\rmdefault}{\mddefault}{\updefault}{ $M_3$}%
}}}}
\put(1951,-5611){\makebox(0,0)[lb]{\smash{{\SetFigFont{12}{14.4}{\rmdefault}{\mddefault}{\updefault}{ $S_{1,3}\times\mathbb C^*$}%
}}}}
\put(3751,-3211){\makebox(0,0)[lb]{\smash{{\SetFigFont{12}{14.4}{\rmdefault}{\mddefault}{\updefault}{ $\left(\mathbb C^*\right)^2$}%
}}}}
\end{picture}%

\end{figure}

\

\

It follows from  chapter 3 in \cite{thesis} that in this limit, the image in the above diagrams of holomorphic curves with bounded energy and genus will converge to piecewise linear graphs. Moreover, holomorphic curves converge locally to a collection of holomorphic maps into the local models which are the fibers in the limiting exploded $\mathbb T$ fibration.

\newpage

\begin{figure}[htb]
 \includegraphics[scale=.8]{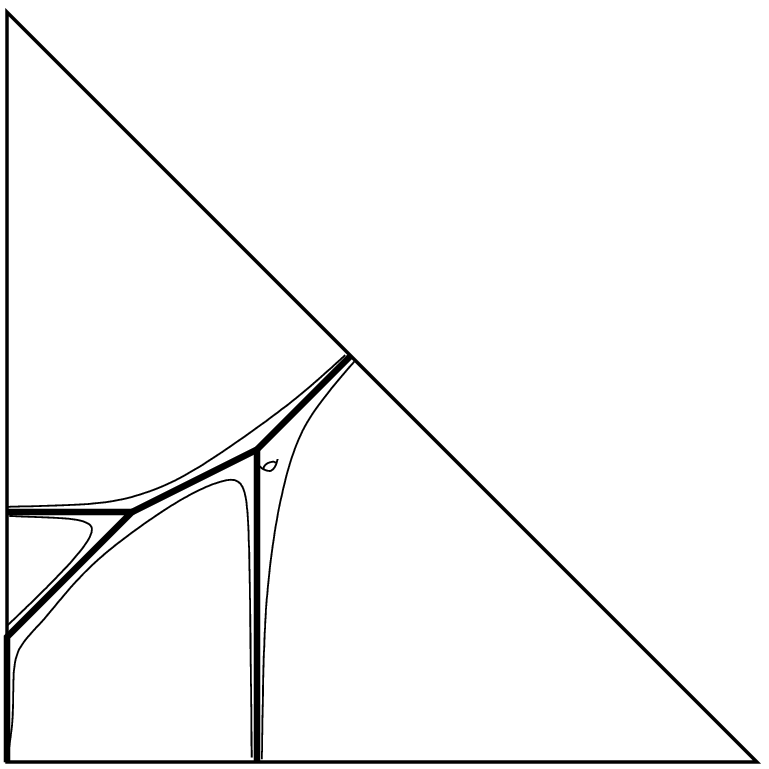}
\end{figure}

The diagram below shows a limiting holomorphic curve, with its image in the base, and schematic pictures of the holomorphic curves in the local model fibers over the base.

\begin{figure}[htb]
\includegraphics[scale=.8]{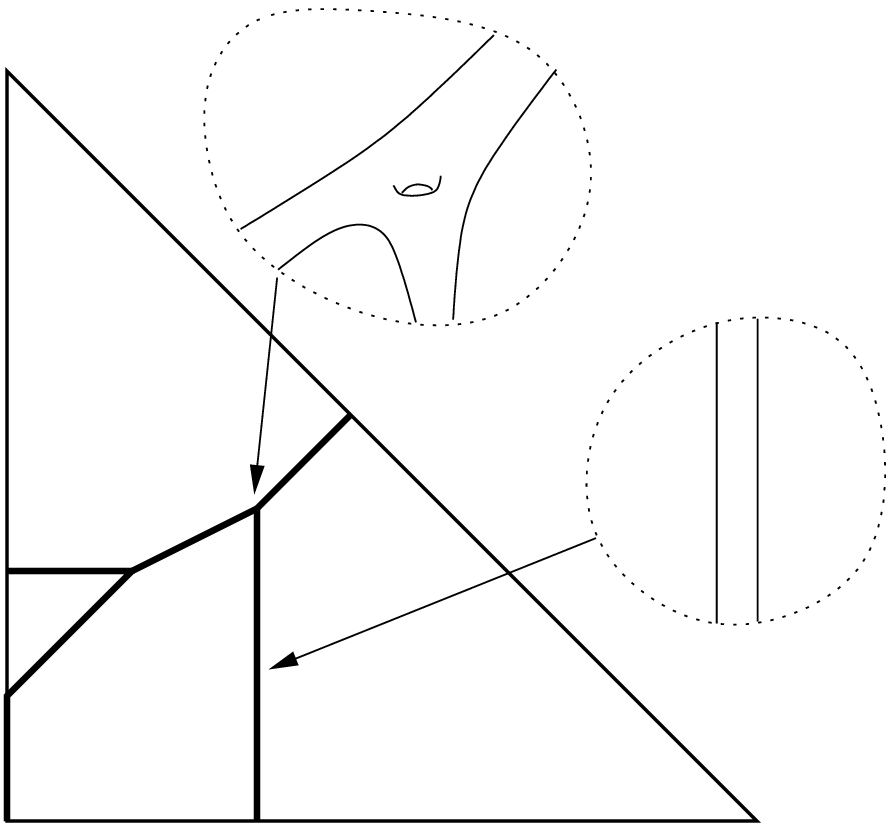} 
\end{figure}

The rest of this paper will be devoted to rigorously defining exploded fibrations and sketching some of their properties. This will not deal with compactness theorems for holomorphic curves and the related questions of how symplectic forms are used to tame holomorphic curves in exploded fibrations.

\

\section{Log smooth structures and the normal neighborhood bundle}\label{log}

The fibers of exploded fibrations will be log smooth spaces.

One way of thinking of a log smooth manifold is a manifold $M$ with boundary and
corners along with the sheaf of vectorfields $C^\infty(\logT M$) which are
tangent to the boundary and corner strata. By manifold with boundary and corners, we mean a manifold with coordinate charts modeled on open subsets of $[0,\infty)^n$. By smooth functions in such a coordinate chart, we mean functions which extend to smooth functions on $\mathbb R^n\supset[0,\infty)^n$.  The sheaf of vectorfields tangent to boundary and corner strata is equal to smooth
sections of a vectorbundle $\logT M$ which we will call the log tangent space. A basis for this vectorbundle in the above coordinate chart is give by $\{x_i\frac\partial{\partial x_i}\}$.
This is often called the $b$ tangent space. The natural objects defined on a manifold with such a structure are smooth sections of tensor products of $\logT M$ and its dual, $\logT^*M$. These tensor fields  have asymptotic symmetry as they approach boundary and corner strata in the sense that the Lie derivative of a smooth tensor field in the direction of $x_i\frac\partial {\partial x_i}$ approaches $0$ like $x_i$. Analysis on manifolds with this structure has been studied extensively by Melrose and others. I wish to thank Pierre Albin for explaining this approach to me. The structure we will describe on log smooth manifolds is the smooth analogue of log structures in algebraic geometry.

The stratified structure of $M$ will be important. Call each connected component of the boundary of $M$ a boundary strata. We will consider this strata to be an immersed submanifold $\iota_S:S\longrightarrow M$ with boundary and corners, which is embedded on the interior of $S$. The other strata of $M$ consist of $M$ itself, and the connected components of higher codimensional corner strata, which are again considered to be immersed submanifolds $\iota_S:S\longrightarrow M$ which are embedded on their interiors.

We will think of log smooth manifolds as possibly non compact manifolds with an
idea of a smooth structure at infinity so that every boundary strata $S$
corresponds to a cylindrical end at infinity $\mathcal N_SM$, called the normal
neighborhood bundle of $S$.

\begin{defn} Given an open subset $U\subset M$ and an open set $U_S\subset S$ inside a boundary strata $\iota_S:S\longrightarrow M$,  a boundary defining function on $U$ defining $U_S$ is a
smooth function $x:M\longrightarrow [0,\infty)$ so that 
\[\{x^{-1}(0)\}=\iota_S(U_S)\subset M\]

and $dx$ is nonvanishing near $S$. 
\end{defn}

 \begin{defn}
 Given an open set $U\subset M$, a log $C^k$ function $f\in\logC^k(U)$  is a function defined on the interior of $U$
   so that 
 \[f:=g+\sum \alpha_i\log(x_i)\]
 where $x_i$ are boundary defining functions on $U$, $\alpha_i\in \mathbb Z$, and $g\in
C^k(U)$.

This defines a sheaf on $M$ of log $C^k$ functions, which we denote by $\logC^k(M)$ 
\end{defn}

\begin{defn}
 A log $C^k$ map between log smooth manifolds
 \[f:M\longrightarrow N\]
 is a continuous map $f$ so that the pull back of the sheaf of log $C^k$
functions on $N$ is contained in the sheaf of log $C^k$ functions on $M$.
\[f^*\left(\logC^k(N)\right)\subset\logC^k(M)\]
\end{defn}

\

 We want to  now enlarge the log smooth category to include integral affine
bundles over log smooth manifolds. 

\begin{defn}
 An integral affine structure on a manifold $M^n$ is a $\mathbb Z^n$ lattice inside   $T_xM$ for each $x\in M$ which is  preserved by some flat connection.
 
  (By an integral $\mathbb Z^n$ lattice, we mean the image of some proper homomorphism $\mathbb Z^n\longrightarrow T_xM$.)

  \
  
  An integral vector in $T_xM$ is a vector in this integral lattice. 
  A basic vector in $T_xM$ is an integral vector which is not a multiple of another integral vector.

\
  
  An integral affine map $f:M^m\longrightarrow N^n$ between integral affine manifolds is a smooth map so that $df$ sends integral vectors to integral vectors.

\end{defn}

\


We will regard a log smooth function as a log smooth map from $M$ to $\mathbb R$
with the standard integral affine structure given by the lattice $\mathbb Z\subset
T\mathbb R$. 
A log smooth map $f:M\longrightarrow \mathbb R^k$ with the standard integral affine structure given by $\mathbb Z^k\subset T\mathbb R^k$ is given by
$(f_1,\dotsc,f_k)$ where $f_i\in {}^\text{log}C^\infty(M)$.

\begin{defn}
 An  integral affine $\mathbb R^k$ bundle $M$ over a log smooth manifold $S$
 is a fibration

\[\begin{split}
   \mathbb R^k\longrightarrow  & M
   \\ & \downarrow\pi
   \\ & S
  \end{split}
\]

so that there exists an open cover $U_\alpha$ of $S$ so that
$\pi^{-1}\left(U_\alpha\right)=\mathbb R^k\times U_\alpha$ and transition
functions 
on $U_\alpha \cap U_\beta$ are of the form
\[(x,s)\mapsto (A(x)+f(s),s)\]
where $f:U_\alpha\cap U_\beta\longrightarrow \mathbb R^k$ is a log smooth map and $A$ is an integral affine automorphism $\mathbb R^k\longrightarrow\mathbb R^k$.
\end{defn}

\

We will call integral affine bundles over log smooth manifolds log smooth
spaces. A strata of a log smooth space is the pullback of this bundle to a strata of the base.

\begin{defn}
 A log $C^k$ function $f\in \logC^k(M)$ on a log smooth space $M$ which is an integral affine $\mathbb R^n$ bundle is a function defined on the interior of $ M$
so that in any local product neighborhood $\mathbb R^n\times U\subset M$
 \[f(x,u):=A(x)+g(u)\]
 where $A:\mathbb R^n\longrightarrow \mathbb R$ is given by an integer matrix
and $g\in\logC^k(U)$. 
\end{defn}

\begin{defn}
 A log $C^k$ map between  log smooth spaces \[f:M\rightarrow N\] 
 is a continuous map from the base of $M$ to the base of $N$ which lifts to a map from the interior of $M$ to the interior of $N$
 so that the pull back of the sheaf of log $C^k$ functions on $N$ is contained in the sheaf of log $C^k$ functions on $M$.
 \[f^*\left(\logC^k( N)\right)\subset\logC^k( M)\]
\end{defn}

\

A log smooth vectorfield on a log smooth space, $v\in C^\infty\left(\logT M\right)$ is a
smooth vectorfield which is tangent to the boundary strata and invariant under
local affine translation of the fibers. These are smooth sections of the log tangent bundle $\logT M$ invariant under local affine translation. Log smooth isotopys of $M$ are generated by the flow of log smooth vectorfields.  

The fact that the flow of log smooth vectorfields preserves boundary strata has
the consequence that we can not naturally identify a neighborhood of a boundary
strata with the log normal bundle. Instead, we have a naturally defined affine
bundle called the normal neighborhood bundle which has an action of the log
normal bundle on it.

\begin{defn}
The normal neighborhood bundle $\mathcal N_sM$ over a point $s\in M$ is the
space of 
functions $\nu_s$ called evaluations
\[\nu_s:{}^\text{log}C^\infty (M)\longrightarrow \mathbb R\]
so that 
\[\nu_s(f)= f(s)\text{ if } f(s)\text { is defined }\]
\[\nu_s(f+g)=\nu_s(f)+\nu_s(g)\]
\end{defn}

\

If $s$ is in a codimension $k$ boundary or corner strata given locally by the vanishing of boundary defining functions $x_1,\dotsc,x_k$, then $\mathcal N_sM$ is an integral affine space equal to $\mathbb R^k$. This identification is given by
$\nu\mapsto\left(\nu(x_1),\dotsc,\nu(x_k)\right)$.

Any log smooth function $f$ induces a function on the normal neighborhood bundle called its restriction by $f(\nu_s):=\nu_s(f)$.

\begin{defn}
The normal neighborhood bundle of a codimension $k$ strata $S\subset M$, $\mathcal N_S(M)$ is a log smooth $\mathbb R^k$ bundle over $S$ with interior given by evaluations $\nu_s\in \mathcal NM$ so that $s$ is in the interior of $S$. The sheaf of log smooth functions  on $\mathcal N_SM$ is given by the restriction of the sheaf of log smooth functions on $M$ in the sense that $f$ is log smooth on $U\subset\mathcal N_SM$ if there exists an open set $\tilde U\subset M$ which contains the projection of $U$ to $S$ and a log smooth function $\tilde f$ defined on $U$ which restricts to $f$.

\end{defn}

\begin{defn}
 An outward direction in $\mathcal N_SM$ is a direction in which the evaluation of any boundary defining function for $S$ decreases.
\end{defn}

This is called `outward' because if we use $N_SM$ to give coordinates on a neighborhood of $S\subset M$,  these directions move out towards the strata $S$ which we consider to be out at infinity.

When all strata $\iota_S:S\longrightarrow M$ are injective, $\mathcal N_SM$ is easier to describe. In this case, $S$ is defined by $k$ boundary defining functions $x_1\dotsc x_k$ and $\mathcal N_SM$  is equal to $\mathbb R^k\times S$ with coordinates $(\nu_1,\dotsc,\nu_k,s)$. The restriction of  $\log(x_i)$ to $\mathcal N_SM$ is given by the coordinate function $\nu_i$. This means that for $g\in \logC^\infty(M)$ defined on the interior of $S$, the restriction of $g(x,s)+\sum \alpha_i\log(x_i)$ is
given by $g(0,s)+\sum\alpha_i\nu_i$. An outward direction in this situation would be any direction that preserves the coordinate $s$ and is negative in the first $k$ coordinates.

One trivial example of a normal neighborhood bundle is $\mathcal N_MM=M$.

\begin{lem}\label{normal restriction}
 Given a log smooth map $f:M\longrightarrow N$, there is a natural map of normal
neighborhood bundles
 \[f: \mathcal N M\longrightarrow \mathcal NN\]
 defined by
 \[f(\nu)(g):=\nu(g\circ f)\]
 Also, if the interior of a strata $S\subset M$ is sent to the interior of a
strata $S'\subset M$, then this is a natural map called the restriction of $f$
to $\mathcal N_SM$ which is a log smooth map
 \[f:\mathcal N_SM\longrightarrow \mathcal N_{S'}N\]
 
\end{lem}

\begin{proof} 

The naturality of this map is clear from the definition. We must check that the induced map  $f:\mathcal N_SM\longrightarrow \mathcal N_{S'}N$ is log smooth. 

First, note that there is a unique continuous map $S\longrightarrow S'$ which is equal to $f$ on the interior of $S$. This gives the continuous map on the base of $\mathcal N_SM$. To see that the pullback of log smooth functions on $\mathcal N_{S'}N$ are log smooth, note that they come locally from the restriction of log smooth functions $g\in \mathcal C^\infty(N)$, so our function is locally $\nu(g)$. Pulling this back via $f$ gives  the function $\nu( g\circ f)$. But $f:M\longrightarrow N$ is log smooth, so $g\circ f$ is log smooth, and therefore the restriction $\nu(g\circ f)$  must also be log smooth.

\end{proof}

In particular, the action on $M$ given by the flow of  a log smooth
vectorfield lifts to an action on $\mathcal NM$. Differentiating this action
gives a lift of log vectorfields $M$ to vectorfields on $\mathcal NM$. This lift
preserves the Lie bracket, addition of vectorfields and multiplication by
functions. It sends nonzero log smooth vectorfields to nonzero vectorfields.

\newpage

 \begin{defn}
  A log smooth morphism \[f:M\rightarrow N\] is a log smooth map 
  \[f:M\rightarrow \mathcal N_SM\] for some strata $S\subset M$. 
 \end{defn}

 \
 
 This finally defines the category of log smooth spaces. We will keep the
distinction between log smooth maps and log smooth morphisms.
 
\section{Stratified integral affine spaces}  \label{base}
  
The base of our exploded fibration will be a stratified integral affine space.

 \begin{defn} A stratified integral affine space is a finite category $B$ associated to a
topological space $\abs B$ with the following structure
 \begin{enumerate}
 
 \item Objects are integral affine strata $B_j^k$ with integral affine
boundary and corners, so that each point in $B_j^k$ has a neighborhood equal to
an open set in $\mathbb R^m\times[0,\infty)^{k-m}$ with the standard integral
affine structure.
 
 \item Morphisms are integral affine inclusions \[\iota:B_i^m\hookrightarrow
B_j^k\] so that each connected affine boundary or corner strata of $B_j^k$ is
the image of a unique inclusion. These inclusions and the identity maps are the
only morphisms.

\item The space $\abs B$ is the disjoint union of strata quotiented out by all
inclusions.
\[\abs B:=\frac{\coprod B_j^k}{\iota: B_i^m\hookrightarrow B_j^k}\]
 \end{enumerate}

 \end{defn}

 \begin{figure}[htb]
\includegraphics[scale=.7]{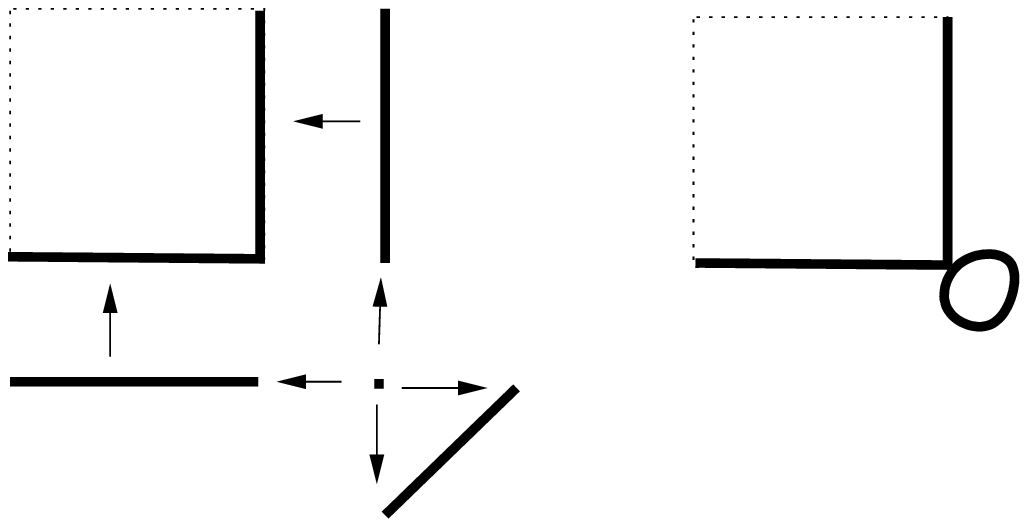}%
%

$B$\hspace{4cm}$\abs B$
\end{figure}

\begin{defn}

A stratified integral affine map from one stratified integral affine space to
another 
\[f:B\longrightarrow C\]
consists of

\begin{enumerate}
\item A functor \[F:B\longrightarrow C\] 

\item Integral affine maps \[f:B_j^k\longrightarrow F(B_j^k)\] so that the
interior of strata in $B$ is sent to the interior of a strata in $C$ and 
 so that \[f\circ \iota:=F(\iota)\circ f\]
\end{enumerate}
\end{defn}

This induces a map on the underlying topological spaces \[f:\abs B
\longrightarrow \abs C\]

 \section{Exploded fibrations}\label{exfib}
 
 \begin{defn}
  An exploded fibration $\ex B$ has the following structure.
  \begin{enumerate}
  \item A stratified integral affine base $B$
  \item A fibration over each strata $\ex B_j^k\longrightarrow B_j^k$ with log smooth fibers \[ M_b:= T_b B_j^k\times  W_b^{n-k}\]  which are integral affine $T_b B_j^k$ bundles. The identification of $M_b$ as a product is not canonical, however the action of $T_bB_j^k$ is well defined. 
  
  \item A flat connection on each fibration $\ex B_j^k$ so that parallel transport gives log smooth isomorphisms compatible with the action of $T B_j^k$ and affine parallel transport.
  
%
 \item For each inclusion $\iota$ in $B$, there is a dual log smooth morphism,
 \[\iota^{\dag}:M_{\iota(b)}\longrightarrow M_b\]
 
 so that 
 \begin{enumerate}
 
 \item For some strata $S\subset M_b$, $\iota^\dag$  is a log smooth isomorphism
  \[\iota^\dag: M_{\iota(b)}\longrightarrow \mathcal N_{ S}
M_b\]

  \item The induced action of $ T_{\iota(b)}B$ on $\mathcal N_{S}M_b$ is
compatible with $d\iota$ in the sense that the induced action of $d\iota(v)$ is equal to the action of $v$. The action of a vector in $T_{\iota(b)}B$ which points towards the interior of the strata  moves $\mathcal N_SM_b$ outwards.
   
 \item Every strata of $M_b$ is the image of a unique isomorphism $\iota^\dag$.
 
  \item \[(\iota_1\circ\iota_2)^\dag=\iota_2^\dag\circ\iota_1^\dag\]
 
 \item The isomorphisms $\iota^\dag$ are compatible with parallel transport.
 \end{enumerate}
  
   \end{enumerate}
  \end{defn}
  
  Note that this definition assumes that the normal neighborhood bundles $\mathcal N_SM_b$ are trivial. This is a global condition. If it is removed, then the base becomes an orbifold version of a stratified integral affine space and the fibration over each strata is a kind of Seifert fibration. 
 
 \begin{defn}
 A  morphism $f:\ex B\longrightarrow \ex C$ is given by
 \begin{enumerate}
 \item A stratified integral affine map of bases $f:B\longrightarrow C$.
  
  \item Log smooth maps \[f_b:M_b\longrightarrow M_{f(b)}\] 
  \begin{enumerate}
  \item
  equivariant with respect to the action of $T_bB$ and $df$
  \item compatible with
the isomorphisms $\iota^\dag$ in the sense that 
  \[f_b\circ\iota^\dag=F(\iota)^\dag\circ f_{\iota(b)}\]
\item compatible with parallel transport within strata.
\end{enumerate}

 \end{enumerate}

 \end{defn}
 
\begin{defn}
A morphism $f:\ex B'\longrightarrow \ex B$ is a refinement if
\begin{enumerate}
\item the map on bases $f:\abs{B'}\longrightarrow\abs B$ is a homeomorphism and
$df$ is an integral affine isomorphism onto its image
\item the maps of fibers $f_b:M_b\longrightarrow M_{f_b}$ are diffeomorphisms of
the interior of $M_b$ to the interior of $M_{f_b}$.
\end{enumerate}\end{defn}

\

\subsection{Coordinates}

\begin{defn}
A log smooth coordinate chart on a log smooth space $M$ which is an integral affine $\mathbb R^m$ bundle is an open
subset $U\subset M$ and log smooth functions
$x_i\in \logC^\infty(U)$, so that 
\begin{enumerate}
\item $e^{-x_i}$ give smooth coordinate functions for $U$. \label{coord 1}
\item\label{coord 2} $\{dx_1,\dotsc,dx_m\}$ restricted to $\mathbb R^m$ fibers
provide a dual basis for the integral affine structure. Call these the affine
coordinate functions.
\item\label{coord 3} $x_i$ for $i>m$ are $\mathbb R^m$ invariant. 
\item \label{coord 4}$\{e^{-x_{m+1}},\dotsc,e^{-x_{m+k}}\}$ are boundary
defining functions which define the boundaries of $U$. Call $x_i$  boundary
coordinate functions for $m<i\leq m+k$. 
\item\label{coord 5} $x_i$ for $i>m+k$ are called the smooth coordinate
functions.

\end{enumerate}

\end{defn}

Another way to say this is that log smooth coordinates on $U$ are induced by a log smooth map $U\longrightarrow \mathbb R^m\times([0,\infty)^k\times\mathbb R^n)$ which is an isomorphism onto its image.
The coordinate functions of $\mathbb R^m$ give the affine coordinates, the boundary coordinate are given by $-\log$ of the coordinate functions of $[0,\infty)^k$, and the smooth coordinate functions are the coordinates of $\mathbb R^n$.

 \
 
 The following lemma follows directly from the definition of log smooth maps.  
 
 \begin{lem}\label{coordinate smooth}
 A map defined on the interior of a log smooth space is log smooth if and only
if the following hold
 \begin{enumerate}
 \item The pull back of affine coordinate functions are log smooth functions,
which means that they are integral sums of affine and boundary coordinate
functions and a smooth function depending on smooth coordinates and boundary
defining functions. \label{coordinate 1}
 \item The pull back of any boundary coordinate function  consists of a positive
integral sum of boundary coordinate functions, and a smooth function depending
on smooth coordinates and boundary defining functions.\label{coordinate 2}
 \item The pull back of any smooth coordinate function consists of a smooth
function depending on smooth coordinates and boundary defining
functions.\label{coordinate 3}
 \end{enumerate}
 \end{lem}
 
 Given a log smooth coordinate chart over $U$ we can canonically associate an
exploded fibration called the explosion of $U$, $\text{Expl}(U)$. The base of
Expl$(U)$ is
 a cone inside the integral affine space $\mathbb R^{m+k}$ given by the subset
where $y_i\geq 0$ for $m<i\leq m+k$.  Let $I\subset\{m+1,\dots,m+k\}$ be some
subset of the boundary indices  and $S_I$ the strata defined by
$\{e^{-x_i}=0,i\in I\}$. The fibers over the set where $\{y_i>0,i\in I\}$ is
$\mathcal N_{S_ I}U$. We can restrict our log smooth coordinate chart on $U$ to
$\mathcal N_{S_I}U$. This changes the boundary coordinates $x_i$ into affine
coordinates for $i\in I$. The action of a base tangent vector $\sum \alpha_i
\frac\partial{\partial y_i}$ on a fiber is given in these coordinates by
$x_i\mapsto x_i+\alpha_i$. 
  
  The construction of Expl$(U)$ is natural in the sense that given any log smooth map $U_1\longrightarrow U_2$ there is a natural map Expl$(U_1)\longrightarrow$Expl$(U_2)$. This tells us that there is a natural construction of of Expl$(M)$ for
any log smooth manifold $M$.

\begin{defn}
 The exploded fibration $\ex R^n$ is defined as the standard fibration with base $\mathbb R^n$ and fibers all equal to $\mathbb R^n$.
\end{defn}

\begin{defn}
 An exploded coordinate chart $\ex U$ on $\ex B$ is a pair of smooth exploded morphisms
 \[\begin{split}&\ex U\longrightarrow \ex B
\\&\downarrow
\\&\ex R^n
\end{split}\]
 which, when restricted to any fiber of $\ex B$ give a log smooth coordinate chart in a suitable choice of basis for $\ex R^n$.
\end{defn}

\begin{defn}
 An exploded fibration $\ex B$ is covered by exploded coordinate charts $\ex U_\alpha$ if the log smooth coordinate charts  obtained by restricting $\ex U_\alpha$
cover every fiber.
\end{defn}

Note that things defined by patching together exploded coordinate charts may fail to have the global property that the fibers are trivial affine bundles. Otherwise, they will always be exploded fibrations.

\subsection{Refinements and blowups}\label{refinements}

In this section, we will show that any subdivision of the base of an exploded fibration
defines a unique refinement. This will also be used later when we explore the intersection theory of exploded fibrations. The construction will be natural, so we will work
in local coordinate charts on the fibers.
  
\

Given the explosion Expl$(U)$ of some log smooth coordinate chart, we will now construct a refinement Expl$(\tilde U)\longrightarrow$Expl$(U)$ for
any integral subdivision of the base into a union of cones. We will use this as a local
model for refinements.

 A cell of this subdivision will be of the form $\{y\cdot\alpha_j\geq 0\}$ for
some collection of $\alpha_j$ which form an integral basis for some subspace 
$\mathbb R^{k+l}\subset\mathbb R^{m+k}$. Any other cell is given by
$\{y\cdot\alpha'_j\geq 0\}$ where $\{\alpha'_j\}$ form a different basis for the
same subspace $\mathbb R^{k+l}\subset\mathbb R^{m+k}$. Each new cell will
correspond to a new log smooth coordinate chart $U^{\{\alpha_i\}}$ with
coordinates
\[\{x\cdot\beta_1,\dotsc,x\cdot\beta_{m-l},x\cdot\alpha_1,\dotsc,x\cdot\alpha_{
k+l},x_{m+k+1},\dotsc\}\] 
 
Here $\{\beta_i,\alpha_j\}$ form an integral basis for $\mathbb R^{m+k}$.
$x\cdot\beta_i$ are the affine coordinate functions, $x\cdot\alpha_i$ are the
boundary coordinate functions, and the remaining $x_i$ are the smooth coordinate
functions. 

The domain of definition of these coordinates, $U^{\{\alpha_i\}}$, is considered
to extend to boundaries where $e^{-x\cdot\alpha_i}=0$. In essence, the new log
smooth coordinates are the old ones with a relabeling of what is considered to
be a boundary coordinate function or affine coordinate function. This operation
is a blowup which creates new boundary components. Each new edge
$\{\alpha_j=0,\text{ for }j\neq i\}$ of a  cell $\{\alpha_j\cdot x\geq0\}$
corresponds to a new boundary component so that the boundary coordinate function
is given by $\alpha_i\cdot x$ and the remaining coordinate functions give
coordinates on the boundary.

The intersection of two of these new coordinate charts include the interiors of
the charts and the boundaries corresponding to common edges. The restriction of
log smooth coordinates to this intersection changes the role of the other
boundary defining functions to smooth coordinate functions.  A transition map on
this intersection between two sets of coordinates then sends affine coordinate
functions to affine coordinate functions, boundary coordinate functions to
boundary defining functions plus some smooth coordinate functions, and smooth
coordinate functions to integral combinations of smooth coordinate functions. As
such they are log smooth. Patching these coordinate charts together gives a log
smooth blowup $\tilde U\longrightarrow U$ . Due to the naturality of the explosion construction, this gives a refinement Expl$(\tilde U)\longrightarrow$Expl$(U)$. 

\begin{lem} \label{refinement model}
 Given any refinement $f:\ex B'\longrightarrow \ex B$, any log smooth coordinate
chart $U$ on a fiber $M_{f_b}$ of $\ex B$ lifts as above to $\tilde
U=f_b^{-1}(U)$.
   
\end{lem}

\begin{proof} 

We need to show that for any log smooth coordinate chart $U$ on $ M_{f(b)}$,
$U^{\{\alpha_i\}}$ gives a log smooth coordinate chart on $f_b^{-1}(U)\subset
M_b$, where $U^{\{\alpha_i\}}$ is the blowup of coordinate charts given by a
cell of the subdivision of the base in the notation above. Moreover, we need to
show that the collection of sets $U^{\{\alpha_i\}}$ from all cells covers $f_b^{-1}U$.
  
The proof is by induction. This holds for fibers that have no boundary
coordinates, as then $f$ is a diffeomorphism and preserves the integral affine
structure on fibers, so it is a log smooth isomorphism, and
$U^{\{\alpha_i\}}=U$.
 
 Suppose that the lemma holds for  fibers of $\ex B'$ with strata of codimension
less than $k$. 

Let $ M_b$ be a fiber of $\ex B'$ with strata of codimension less than or equal
to $k$, and let $ U$ be a coordinate chart on $M_{f(b)}$. The subdivision of the
base $B$ gives a subdivision of the base of Expl$(U)$. Choose a cell of this
subdivision given by $\{\alpha_i\geq 0\}$. We must show that $U^{\{\alpha_i\}}$
gives a log smooth coordinate chart on $M_b$.

Let $I$ be a proper subset of the indexes of $\{\alpha_i\}$,  $S_I\subset M_b$
and $\iota_I$ the corresponding strata and base inclusion. Denote the
restriction of $U$ to $\mathcal N_{f(S_I)}M_{f(b)}$ by 
$U_I$.
 $\mathcal N_{S_I}M_b=M_{\iota_I(b)}$ only has strata of codimension less than
$k$ so the inductive hypothesis holds and $U_I^{\{\alpha_i,i\in I\}}$ is a log
smooth coordinate chart on $\mathcal N_{S_I}M_b=M_{\iota_I(b)}$. This chart
$U_I^{\{\alpha_i,i\in I\}}$ is also the restriction of $U^{\{\alpha_i\}}$
 to $\mathcal N_{S_I}M_b$, so we have that the restriction of $U^{\{\alpha_i\}}$
to the normal neighborhood bundle of any boundary strata gives log smooth coordinates. We also have that the
restriction of all such $U^{\{\alpha_i\}}$ to normal neighborhood bundles
covers the restriction of $f_b^{-1}(U)$ to the normal neighborhood bundle.
Therefore  the sets $U^{\{\alpha_i\}}$ must cover $f^{-1}(U)$ because each one
covers the interior, and together they cover the boundary strata. 

Let's now check that the axioms of a log smooth coordinate chart are satisfied by
$U^{\{\alpha\}}$. First, all coordinates are log smooth functions because $f$ is
a log smooth map, and the coordinates are just pullbacks of integral sums of log
smooth functions.

Item \ref{coord 2} is satisfied, as our choice of fiber coordinate functions is
just the pullback of a subset of the fiber coordinate functions on the target,
$f_b$ is equivariant with respect to the action of the tangent space of the base
and $df$, and $df$ is an integral affine isomorphism onto its image.

Item \ref{coord 3} is satisfied.  The boundary coordinates $x\cdot\alpha_i$ are
constant on fibers  because $f_b$ is equivariant and the $\alpha_i$ are
orthogonal to the fiber directions. The smooth coordinates are just pullbacks of
smooth coordinates, so they too are constant on fibers.

Item \ref{coord 4} is satisfied. $e^{-x\cdot\alpha_i}$ is a boundary defining
function because it is log smooth, and restricts correctly to the normal
neighborhood of the boundary it defines. 

Item \ref{coord 1} is satisfied. This is because the coordinate functions are
log smooth, give smooth coordinates on the interior of $U^{\{\alpha_i\}}$, their
restrictions to normal neighborhood bundles give smooth coordinate functions,
and $e^{-x\cdot\alpha_i}$ are boundary defining functions. 

Item \ref{coord 5} is an empty condition, so $U^{\{\alpha_i\}}$ is a log smooth
coordinate chart and the lemma is proved.

\end{proof}

\begin{prop}\label{lift}
Given a  morphism $f:\ex B\rightarrow \ex C$ and a refinement  $\ex
C'\longrightarrow \ex C$ so that the map of bases $f:B\rightarrow C$ lifts to a
map
$\tilde f:B\longrightarrow C'$, there exists a unique lift  $\tilde f:\ex
B\longrightarrow \ex C'$ of $f$ to a morphism to $\ex C'$.  
\end{prop}

\begin{proof} 

Refinements are diffeomorphisms on the interiors of fibers, and log smooth maps
are determined by their restriction to the interior, so uniqueness is automatic.
We just need to check that the maps given by $f_b=\tilde f_b$ restricted to the
interior of $M_b$ are log smooth and satisfy the compatibility requirements for
exploded fibrations.

To show that $\tilde f_b$ is log smooth, we use the criteria from Lemma
\ref{coordinate smooth} and the local normal form for refinements from Lemma
\ref{refinement model}. 

Affine coordinates on the blowup $\tilde M_c$  are a subset of affine coordinate
on $M_c$, so criterion \ref{coordinate 1} is automatically satisfied. 
 Smooth coordinates for $\tilde M_c$ in the normal form from Lemma
\ref{refinement model} are just smooth coordinates from $M_c$, so criterion
\ref{coordinate 3} is also automatically satisfied, so we just need to study the
pull back of boundary coordinate functions.
 
The fact that $f$ is log smooth tells us the following: the pull back under
$\tilde f$ of boundary coordinate functions on $\tilde M_c$ are integral sums of
affine and boundary coordinates plus a smooth function depending on boundary
defining functions and smooth coordinates. The exact combination of affine and
boundary coordinates is determined by the map $f:B\longrightarrow C$ on bases. The
requirement that $f$ lifts to the subdivision $C'$ is equivalent to this
combination being a positive combination of boundary coordinates, which tells us
that 
$\tilde f$ satisfies criterion \ref{coordinate 3}, and is therefore log smooth.

 The fact that $\tilde f$ is compatible with parallel transport follows
from the fact that $ f$ is compatible and $\tilde f$ and $f$ are equal on the
interiors of fibers. Also, restricted to interiors, $\tilde f\circ
\iota^\dag=f\circ \iota^\dag=F(\iota)^\dag\circ f= \tilde F(\iota)^\dag\circ
\tilde f$, so $\tilde f$ is also compatible with the inclusions $\iota^\dag$ and
$f$ is an exploded morphism. 
 
 \end{proof}

\begin{thm}\label{subdivision}
Given an exploded fibration $\ex B$ and a stratified integral affine map
\[f:B'\longrightarrow B\]
so that $f$ is a homeomorphism and $df$ is an integral affine isomorphism onto
its image, there exists a unique refinement $ f:\ex B'\longrightarrow \ex B$
with base $f:B'\rightarrow B$.  
\end{thm}

\begin{proof} 

First, Proposition \ref{lift} tells us that this refinement is unique up to
isomorphism.

To show the existence of a refinement, note that any exploded fibration is
locally modelled on  Expl$(U)$ for some log smooth coordinate chart $U$. Given a
subdivision, we constructed a refinement  Expl$(\tilde
U)\longrightarrow$Expl$(U)$. This is a local model for our refinement. These
refinements so constructed must coincide on the intersection of their domains of
definition due to the uniqueness of refinements from Proposition \ref{lift}, so
these local refinements patch together to a global refinement.

\end{proof}

 \section{Exploded $\mathbb T$ fibrations}\label{etf}

Exploded $\mathbb T$ fibrations can be thought of as exploded fibrations with
extra torus symmetry so that instead of fibers being $\mathbb R^k$ bundles,
fibers are $(\mathbb C^*)^k$ bundles. We will first describe the analogue of log
smooth spaces in this setting.
 
 \subsection{Log smooth $\mathbb T$ spaces}

A log smooth $\mathbb T$ manifold is a connected smooth manifold $M$ locally modeled on $\mathbb C^k\times \mathbb R^{n-2k}$ with the following sheaf of log smooth $\mathbb T$ functions.

\begin{defn}
 A log smooth $\mathbb T$ function $f\in\TC^\infty(\mathbb C^k\times \mathbb R^{n-2k})$ is a map \[f:\left(\mathbb C^*\right)^k\times \mathbb R^{n-2k}\longrightarrow \mathbb C^*\]
   \[f:=z_1^{\alpha_1}\dotsb z_k^{\alpha_k}g\]
   \[\text{for }\alpha_i\in\mathbb Z,\text{ and }g\in C^\infty(\mathbb C^k\times \mathbb R^{n-2k},\mathbb C^*)\]
   
\end{defn}

\begin{defn}
 A log smooth $\mathbb T$ manifold $M$ is a smooth manifold $M$ along with a sheaf
 of log smooth $\mathbb T$ functions $\TC^\infty(M)$, where if $U\subset M$ is an open subset, $f\in TC^\infty(U)$ is some $\mathbb C^*$ valued function defined on a dense open subset of $U$.
   
   This must satisfy the condition that around any point in $M$, there exists a neighborhood $U$ and a diffeomorphism $\phi:U\longrightarrow \mathbb C^k\times \mathbb R^{n-2k}$ so that 
 \[\phi^*\left(\TC^\infty(\mathbb C^k\times\mathbb R^{n-2k})\right)=\TC^\infty(U)\]
\end{defn}

\begin{defn}
 A smooth function $z:U\subset M\longrightarrow\mathbb C$ is a boundary defining function if $dz$ is nonzero when $z$ is zero, and $z\in\TC^\infty(U)$ 
\end{defn}

\begin{defn}
A strata $S\longrightarrow M$ is a connected properly immersed  submanifold of $M$ which is given locally by the vanishing of some number of boundary defining functions, and which is embedded on a dense subset of $S$.

The interior of a strata $S$ is the strata minus all substrata.
\end{defn}

Note that each point in $M$ is in the interior of a unique strata.

\

With log $\mathbb T$ functions defined, other definitions are analogous to the
definitions in the log smooth case.

\begin{defn}
 \[f: M\rightarrow  N\]
 is a log $C^k$ $\mathbb T$ map if it is continuous, sends the interior of  $M$
to the interior of $N$, and 
 \[f^*\left(\TC^k( N)\right)\subset\TC^k( M)\] 
\end{defn}

\begin{defn}
\item A log smooth $\mathbb T$ space is a $\left(\mathbb C^*\right)^n$
bundle
 \[\begin{split}(\mathbb C^*)^n\longrightarrow   &M\\&\downarrow \pi
 \\& S
 \end{split}\]
over a log smooth $\mathbb T$ manifold $ S$.

In particular, there exists an open cover $\{U_\alpha\}$ of $ S$ so that
$\pi^{-1}U_\alpha=(\mathbb C^*)^n\times U_\alpha $ and transition functions are
of the form \[(z,s)\mapsto (\phi(z)f'(s),s)\] where\item $f':U_\alpha\cap U_\beta\longrightarrow
(\mathbb C^*)^n$ is a product of log smooth $\mathbb T$ functions,\item and $\phi:\left(\mathbb C^*\right)^n\longrightarrow \left(\mathbb C^*\right)^n$ is a group isomorphism.

\end{defn}

\

 A $\mathbb C^*$ valued function $f$ is in $\TC^k(\mathbb C^*\times S)$ if 
 $f=z^mg$ where $g\in\TC^k( S)$ and $m\in \mathbb N$. Using this inductively
defines the sheaf of log smooth $\mathbb T$ functions on $(\mathbb C^*)^n\times 
S$, and then using this locally defines the sheaf of log smooth $\mathbb T$
functions
 on a log smooth $\mathbb T$ space $M$ which is a nontrivial $\left(\mathbb C^*\right)^n$ bundle.
 
 \begin{defn}
  \[f:   M\longrightarrow  N\]
  is a log $C^k$ $\mathbb T$ map between log smooth $\mathbb T$ spaces if it
projects to a log $C^k$ $\mathbb T$ map on the bases of $M$ and $N$ , and
  \[f^*\left(\TC^k(  N)\right)\subset \TC^k( M) \]
 \end{defn}

\begin{defn}
 The normal neighborhood bundle over a point $s\in M$, $\mathcal N_sM$ is given
by the space of evaluations $\nu_s:\TC(M)\longrightarrow \mathbb C^*$ so that
 \item \[\nu_s(f)=f(s)\text{ if well defined }\]
 \item \[ \nu_s(fg)=\nu_s(f)\nu_s(g)\] 
 \end{defn}

 \begin{defn}
  The normal neighborhood bundle of a codimension $2k$ strata $S\longrightarrow M$ is a log smooth $\mathbb T$ space which is a $(\mathbb C^*)^k$ bundle over $S$. The fiber in $\mathcal N_SM$ over a point $s$ in the interior of $S$  correspond to evaluations $\nu_s$ over $s$.

If  $S$ is locally given by vanishing of boundary defining
functions $z_1, \dotsc,z_k$, then local coordinates for  $\mathcal N_SM$ are
given by $(\nu(z_1),\dotsc,\nu(z_k), s)$, where the first $k$ coordinates give
coordinates for the $\left(\mathbb C^*\right)^k$ fiber, and $s$ denotes the complimentary coordinates which give coordinates on $S$. 

 \end{defn}

 \begin{defn}
  An outward direction in $\mathcal N_SM_b$ is a direction in which the restriction of any boundary defining function for $S$ decreases. 
 \end{defn}

\begin{lem}
 A log $\mathbb T$ map $f:M\longrightarrow N$ induces a natural map
 \[f:\mathcal N M\longrightarrow \mathcal N N\]
\[\text{defined by } f(\nu)(g)=\nu(f\circ g)\]
If $f$ sends the interior of a strata $S$ to the interior of $S'$, then this
gives a log  $\mathbb T$ map
\[f:\mathcal N_SM\longrightarrow \mathcal N_{S'}N\]
\end{lem}

The proof of this lemma is analogous to Lemma \ref{normal restriction}

\begin{defn}
 A log $\mathbb T$ morphism from $M$ to $N$ is a log $\mathbb T$ map 
 \[f:M\longrightarrow \mathcal N_SN\]
 for some strata $S\in N$.
\end{defn}

If we take the oriented real blowup of all strata of a log smooth $\mathbb T$
space, we obtain a log smooth space. Locally, this involves replacing any torus
boundary defining function $z$ with a log smooth boundary defining function
$\abs z$ and a smooth coordinate $\frac z{\abs z}$. In this way, we can view any
log smooth $\mathbb T$ space as a log smooth space. Any log $\mathbb T$ morphism
is still a log morphism, and the normal neighborhood bundle in the log $\mathbb
T$ setting is the normal neighborhood in the log smooth setting. This allows us
to use any concept from the log smooth setting.

\subsection{Exploded $\mathbb T$ fibrations}

We can now define exploded $\mathbb T$ fibrations analogously to exploded
fibrations.

 \begin{defn}
  An exploded $\mathbb T$ fibration $\ex B$ has the following structure.
  \begin{enumerate}
  \item A stratified integral affine base $B$
  \item For each strata $B_j^k$, a fibration $\ex B_j^k\longrightarrow B_j^k$ with log smooth $\mathbb T$ fibers for each $b\in B_j^k$,
   \[ M_b:= \left(\mathbb C^*\right)^k\rtimes  W_b^{n-2k}\]  These fibers are principal $\left(\mathbb C^*\right)^k$ bundles. They have an action
of $T_b B_j^k$ given in integral coordinates by $(x_1,\dotsc,x_k)$ acts by
multiplication by $(e^{x_1},\dotsc,e^{x_k})$.

  \item A flat connection on each of these fibrations $\ex B_j^k$ so that parallel transport gives log $\mathbb T$ smooth isomorphisms which are compatible with the action of $T_b B_j^k$ and affine parallel transport.
 \item For each inclusion $\iota$ in $B$, there is a dual log smooth $\mathbb T$
morphism,
 \[\iota^{\dag}:M_{\iota(b)}\longrightarrow M_b\]
 
 so that 
 \begin{enumerate}
 
 \item For some strata $S\subset M_b$, $\iota^\dag$  is a  log smooth
$\mathbb T$ isomorphism
  \[\iota^\dag: M_{\iota(b)}\longrightarrow \mathcal N_{ S}
M_b\]

  \item The induced action of $T_{\iota(b)}B$ on $\mathcal N_SM_b$ is
compatible with $d\iota$. The action of any vector in $T_{\iota(b)}B$ pointing towards the interior of the strata moves $\mathcal N_SM_b$ in an outward direction.
   
 \item Every strata of $M_b$ is the image of a unique isomorphism $\iota^\dag$.
 
  \item \[(\iota_1\circ\iota_2)^\dag=\iota_2^\dag\circ\iota_1^\dag\]
 
 \item These isomorphisms $\iota^\dag$ are compatible with parallel transport within strata.
 
 \end{enumerate}
  
   \end{enumerate}
  \end{defn}

  \
   
  As with our definition of exploded fibrations, there is a global assumption that fibers have a globally defined action of $\left(\mathbb C^*\right)^k$ on them. If we patch together exploded $\mathbb T$  coordinate charts, then we will obtain a more sophisticated version of an exploded $\mathbb T$ fibration with an orbifold base and Seifert fibrations over strata.

 \begin{defn}
 A  morphism of exploded $\mathbb T $ fibrations, $f:\ex B\longrightarrow \ex C$ is given by
 \begin{enumerate}
 \item A stratified integral affine map of bases $f:B\longrightarrow C$.
    \item Log smooth $\mathbb T$ maps \[f_b:M_b\longrightarrow M_{f(b)}\] 
   \begin{enumerate}
  \item
  equivariant with respect to the action of $T_bB$ and $df$
  \item compatible with
the isomorphisms $\iota^\dag$ in the sense that 
  \[f_b\circ\iota^\dag=F(\iota)^\dag\circ f_{\iota(b)}\]
\item compatible with parallel transport within strata.
\end{enumerate}
  
 \end{enumerate}

 \end{defn}

 \
 
Note that any exploded $\mathbb T$ fibration is an exploded fibration and that
any exploded $\mathbb T$ morphism is also a morphism of exploded fibrations. 

\

Refinements are defined the same way
as in the log smooth case. The results of section \ref{refinements} hold for exploded $\mathbb T$ fibrations.

\subsection{Almost complex structure}

The tangent space $T\ex B$ of an exploded fibration is an exploded fibration
with the same base and fibers given by $\logT M_b$, the log tangent space of
fibers of $\ex B$.
This has a natural exploded structure because there is a natural identification
of  $\mathcal N_S\logT M_b$ with $\logT \mathcal N_SM_b$. The cotangent space
$T^*\ex B$ is a similarly well defined exploded fibration. A smooth section of
one of these bundles is a section $\ex B\longrightarrow T\ex B$ which is a
smooth exploded morphism. This allows us to define things like metrics or
complex structures.

\begin{defn}
An almost complex structure on an exploded fibration $\ex B$ is an endomorphism
$J$ of $T\ex B$ which squares to minus the identity, given by a smooth section
of $T\ex B\otimes T^*\ex B$. 
\end{defn}

(The tensor product above is over smooth functions on $\ex B$ which are smooth
morphisms to the exploded fibration which has a single point as its base and
$\mathbb R$ as the fiber over that point.)

\begin{defn}
 An (almost) complex exploded $\mathbb T$ fibration $\ex B$ is an almost complex
structure on $\ex B$ considered as an exploded fibration which comes from giving
each fiber $M_b$ the structure of a (almost) complex $(\mathbb C^*)^k$ bundle over an (almost)
complex manifold.
   
\end{defn}

Note that this implies that substrata are holomorphic submanifolds.

\

There should be a good theory of holomorphic curves in almost complex exploded
fibrations. In general, we need to consider holomorphic curves which are morphisms of a weaker type than described above. The theory of holomorphic curve in exploded $\mathbb T$ fibrations is better behaved. In this case, holomorphic curves are exploded $\mathbb T$ morphisms and when transversality conditions are met, the moduli space has an exploded $\mathbb T$ structure.

\section{Examples of exploded $\mathbb T$ curves}\label{examples}

\subsection{Moduli space of stable \exploded curves}

\begin{defn}
An \exploded curve is a complex \exploded fibration with a one dimensional base
so that
\begin{enumerate}
\item The base is  complete when given the metric that gives basic integral
vectors unit size.

\item  The fiber over any vertex is a compact Riemann surface with strata
corresponding to a collection of marked points. 
\end{enumerate}
\end{defn}
 
We will also call a holomorphic morphism of an \exploded curve to
a almost complex \exploded fibration an \exploded curve.

\begin{defn}
 An \exploded curve is stable if it has a finite number of
automorphisms, and it is not a nontrivial refinement of another \exploded curve. 
\end{defn}

\

We now describe in detail the moduli space of stable \exploded curves (mapping
to a point.)

\

Consider a Riemann surface $\Sigma_c$  which is the fiber over the point $c$ in
the base $C$ of an \exploded curve $\ex C$.  This has a number of marked points
which correspond to boundary strata.

If $x\in \Sigma_c$ is such a marked point, then we can identify $\mathcal
N_x\Sigma_c$ with $T_x\Sigma_c-0$. The identification is given as follows: If
$x$ is locally defined by the vanishing of a boundary defining function $z$, a
nonzero vector $v\in T_x\Sigma_c$ corresponds to $\nu\in \mathcal N_x\Sigma_c$
where $\nu(z):=dz(v)$. (This determines $\nu\in\mathcal N_x\Sigma_c$ by
$\nu(z^kg)=(dz(v))^kg(x)$.)

\begin{lem}
A (stable) \exploded curve is equivalent to a (stable) punctured nodal Riemann
surface with the following extra information at each node:

\begin{enumerate}
\item A length $l\in(0,\infty)$

\item A node consists of two marked points which are considered to be joined.
There is a nondegenerate $\mathbb C$ bilinear pairing of the tangent spaces at
these points.

\end{enumerate}

\end{lem}

The length is the length of the edge joining two marked points. The $\mathbb C$
bilinear pairing comes from the identification of normal neighborhood bundles
at these points.

\

\begin{lem}\label{automorphism}
Any finite order automorphism $\phi:\ex C\longrightarrow\ex C$  of a connected exploded $\mathbb T$ curve that is trivial on some fiber is the identity.
\end{lem}

\begin{proof} 

Suppose that $\phi:\ex C\longrightarrow \ex C$ is such an automorphism. If $\phi_v$ is the identity for some vertex $v$, then the restriction of $\phi_v$ to the normal neighborhood bundle is the identity, so $\phi_e$ is also the identity for $e$ any point on an edge connected to $v$ (and $\phi$ is the identity on those edges.)

If $\phi_e$ is the identity for $e$ some point in an edge, $\phi$ must be the identity on the edge, and all fibers on the edge. Suppose that $v$ is a vertex on the end of such an edge. $\phi_v:\Sigma_v\longrightarrow \Sigma_v$ must be a finite order automorphism of $\Sigma_v$ that fixes the marked point $x$ attached to the edge, and that fixes the tangent space $T_x\Sigma_v$. The only such automorphism is the identity. (This can be seen locally, if $\phi_v(z)=z+cz^n+O(z^{n+1})$, then the $k$-fold composition of $\phi_v$ must be $z+kcz^n+O(z^{n+1})$, so if $\phi_v$ is not the identity, then it is not of finite order.) 

\end{proof}

In particular, the above lemma tells us that the automorphisms of nodal Riemann surfaces that are unbalanced in the sense that they twist one side of a node more than the other are not automorphims of exploded $\mathbb T$ curves. An automorphim such as this lifts to an isomorphism of an exploded curve with an exploded curve obtained by modifying the bilinear pairings at nodes. Isomorphims such as this will act like automorphims when we describe the orbifold structure of moduli spaces of exploded curves.

\begin{defn}
 A lifted automorphism of an \exploded  curve is an isomorphism of the \exploded curve to an exploded $\mathbb T$ curve obtained by modifying bilinear pairings at nodes.
\end{defn}

The lifted automorphisms of an \exploded curve are in one to one correspondence with a subgroup of the automorphisms of the underlying nodal curve. In the case that all edges have equal length, the lifted automorphisms correspond to the automorphisms of the underlying nodal curve.

\

The moduli space $\ex M_{g,n}$ of stable \exploded curves with genus $g$ and $n$ marked points has the structure of an orbifold exploded $\mathbb T$ fibration. An orbifold exploded fibration has uniformising charts $(\ex U,G)$ where $\ex U$ is an exploded $\mathbb T$ coordinate chart, and $G$ is a group acting on $\ex U$.
 
 Uniformising charts for $\ex M_{g,n}$ are constructed as follows: Choose holomorphic uniformising charts $(U,G)$ for Deligne Mumford space $\bar{\mathcal M}_{g,n}$ so that the boundary strata contained in $U$ consist of the vanishing of  $\{z_1,\dotsc,z_k\}$, the first $k$ coordinates of $U$. This makes $U$ a log smooth $\mathbb T$ coordinate chart where $\{z_1,\dotsc,z_k\}$ are boundary defining functions. These correspond to $k$ nodes in the curve over $z_1=\dotsb=z_k=0$. Note that the action of $G$ on $U$ and any transitions between coordinate charts are log $\mathbb T$ smooth.
  
  Uniformising charts for $\ex M_{g,n}$ are then given by the explosion $(\text{Expl}(U),G)$. The base of Expl$(U)$ has coordinates  $(l_1,\dotsc,l_k)\in[0,\infty)^k$.  The fiber over the strata where $l_i=0$ for $i\notin I$ is given by $\mathcal N_{S_I}U$ where $S_I$ is the strata given by $z_i=0$ for $i\in I$.
   The action of $G$ and transition maps are given by noting that any log smooth $\mathbb T$ map $U_1\longrightarrow U_2$ gives a natural exploded $\mathbb T$ morphism Expl$(U_1)\longrightarrow$Expl$(U_2)$.  
   
   \
   
   We must now see why points in $\ex M_{g,n}$ (which we will consider to be the image of a morphism of a point to $\ex M_{g,n}$) correspond to exploded $\mathbb T$ curves. We do this for each uniformising chart Expl$(U)$.
   A point in the strata  where  $l_i=0$ for $i\notin I$ has the information of a nodal curve corresponding to $s\in S_I$ plus an evaluation $\nu_s\in \mathcal N_sU$ and lengths $l_i$ for $i\in I$. These lengths correspond to the length associated to each node. Note that if this point is the image of a morphism, then $l_i>0$ for $i\in I$, and these are the only nodes that our curve has. To describe this as a moduli space of exploded curves, we must also say how the bilinear pairing of tangent spaces at a node is given by the above choice of an evaluation  $\nu_s\in\mathcal N_sU$.  
    
    Recall that there is a universal curve $\bar{\mathcal M}_{g,n+1}\longrightarrow \bar{\mathcal M}_{g,n}$ given by forgetting the point labeled with $(n+1)$. This is log $\mathbb T$ smooth as it is holomorphic and the inverse image of boundary strata are boundary strata, so we also have a map of normal neighborhood bundles $\mathcal N\bar{\mathcal M}_{g,n+1}\longrightarrow \mathcal N\bar{\mathcal M}_{g,n}$. If the forgotten point is near the first node, this map can be given in coordinates locally near $z_1=0$ by \[(z_1^+,z_1^-,z_2,\dots)\mapsto (z_1^+z_1^-,z_2,\dots)\] (where dots are completed with the identity). The curve over the point $z_1=0$ gives local coordinates for our node, considered as two disks with coordinates $z_1^+$ and $z_1^-$ joined over the points $z_1^+=0$ and $z_1^-=0$. We have a map of the normal neighborhood bundle of the strata where $z_1^+=z_1^-=0$ to the normal neighborhood bundle of the strata $z_1=0$, given in our coordinates by the above map, where $z_1^+$ and $z_1^-$ now give coordinates for the tangent space to our curve at the node. The bilinear pairing in these coordinates is \[(z_1^+,z_1^-)\mapsto \frac{z_1^+z_1^-}{\nu_s(z_1)}\]
   This is independent of coordinate choices, and the choice of  pairings for each of the $k$ nodes is equivalent to the choice of $\nu_s$.
   
 A second way to see the relationship between pairings and $\nu_s$ is to view the pairings as gluing information. In particular, given the curve over $s$ and choice of coordinates $z_i^+,z_i^-$ around its nodes, we can glue these neighborhoods with the help of a small complex parameter $c_i$, identifying  \[z_i^+= \frac {c_i}{z_i^-}\]  
 This defines a gluing map $\mathbb C^k\longrightarrow U$ defined near $c_i=0$. Giving $\mathbb C^k$ the strata $c_i=0$, this map is log $\mathbb T$ smooth, and the restriction $\mathcal N_0\mathbb C^k\longrightarrow\mathcal N_sU$ is an isomorphism. If the pullback of $\nu_s$ to $\mathcal N_0\mathbb C^k$ has coordinates $(\nu_{s,1},\dotsc,\nu_{s,k})$, then the pairing between $z_i^+$ and $z_i^-$ is given by \[(z_i^+,z_i^-)=\frac{z_i^+z_i^-}{\nu_{s,i}}\]

 The above discussion  can be summarized in the following lemma.
 
 \begin{lem}
  The forgetfull map 
  \[\ex M_{g,n+1}\longrightarrow \ex M_{g,n}\]
 is log $\mathbb T$ smooth. Given any point $p\longrightarrow \ex M_{g,n}$, the fiber over $p$ is the exploded curve corresponding to $p$ quotiented out by its automorphisms. 
 \end{lem}

For example, $\ex M_{0,3}$ is a point, and $\ex M_{0,4}$ is equal to an exploded $\mathbb T$ curve with $3$ punctures.  A curve in $\ex M_{0,4}$ might look like the following curve or a similar one where the punctures 1 and 3 or 1 and 4 are grouped together.

\begin{figure}[htb]
 \begin{picture}(0,0)%
\includegraphics{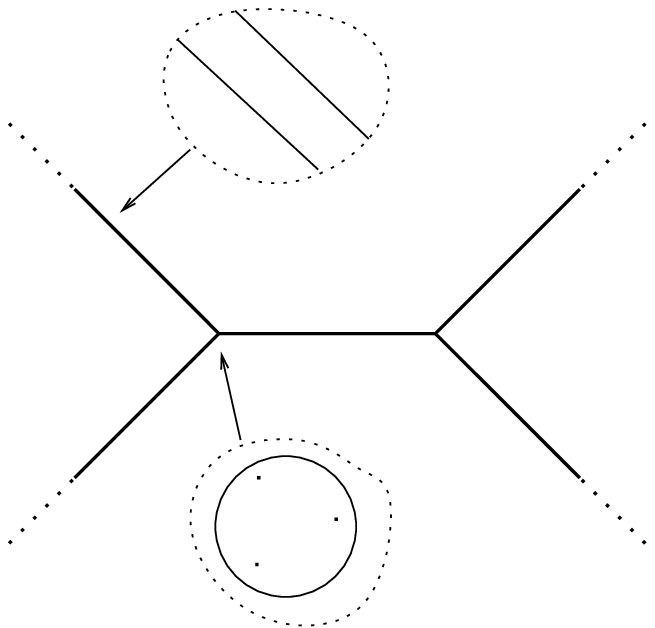}%
\end{picture}%
\setlength{\unitlength}{3947sp}%
\begingroup\makeatletter\ifx\SetFigFont\undefined%
\gdef\SetFigFont#1#2#3#4#5{%
  \reset@font\fontsize{#1}{#2pt}%
  \fontfamily{#3}\fontseries{#4}\fontshape{#5}%
  \selectfont}%
\fi\endgroup%
\begin{picture}(3312,2983)(2973,-3001)
\put(4934,-2574){\makebox(0,0)[lb]{\smash{{\SetFigFont{7}{8.4}{\rmdefault}{\mddefault}{\updefault}{ $\mathbb CP^1-\{0,1,\infty\}$}%
}}}}
\put(4942,-263){\makebox(0,0)[lb]{\smash{{\SetFigFont{7}{8.4}{\rmdefault}{\mddefault}{\updefault}{ $\mathbb C^*$}%
}}}}
\put(3090,-471){\makebox(0,0)[lb]{\smash{{\SetFigFont{7}{8.4}{\rmdefault}{\mddefault}{\updefault}{ 1}%
}}}}
\put(3142,-2704){\makebox(0,0)[lb]{\smash{{\SetFigFont{7}{8.4}{\rmdefault}{\mddefault}{\updefault}{ 2}%
}}}}
\put(6232,-653){\makebox(0,0)[lb]{\smash{{\SetFigFont{7}{8.4}{\rmdefault}{\mddefault}{\updefault}{ 3}%
}}}}
\put(6016,-2384){\makebox(0,0)[lb]{\smash{{\SetFigFont{7}{8.4}{\rmdefault}{\mddefault}{\updefault}{ 4}%
}}}}
\put(4321,-1351){\makebox(0,0)[lb]{\smash{{\SetFigFont{12}{14.4}{\rmdefault}{\mddefault}{\updefault}{ $l$}%
}}}}
\end{picture}%

\end{figure}

For each picture of the above type, there is a $\mathbb C^*$ family of curves in $\ex M_{0,4}$ corresponding to different identifications of normal neighborhood bundles along the interior edge. The length $l$ of this edge gives a coordinate for a strata of $M_{0,4}$. The above observation tells us that the fiber over each point in this strata is $\mathbb C^*$. The different pairings of punctures give $3$ such strata. the end $l=0$ of this strata is glued to a strata with curves looking like the following:

\newpage

\begin{figure}[htb]
 \includegraphics[scale=.6]{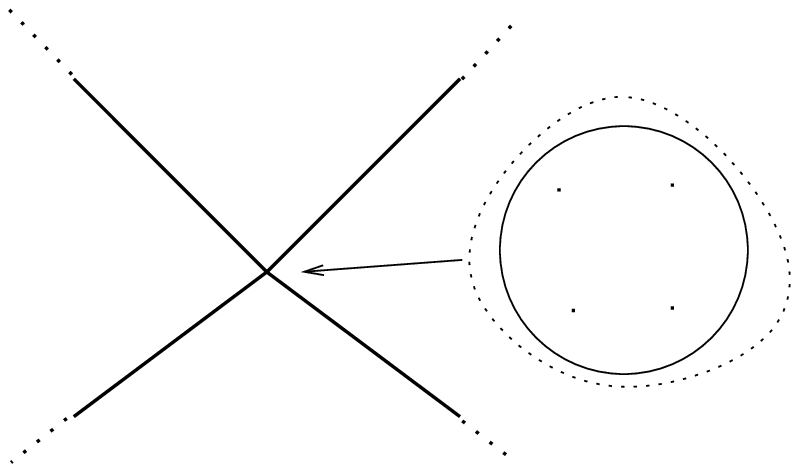}
\end{figure}

The curves with a picture such as this form a family parametrized by the complex structure of the fiber over the central vertex. This gives a strata of $\ex M_{0,4}$ consisting of a point with fiber equal to $\bar {\mathcal M}_{0,4}$ with $3$ boundary strata. A picture of $\ex M_{0,4}$ is as follows
 
\begin{figure}[htb]
 \begin{picture}(-4,0)%
\includegraphics{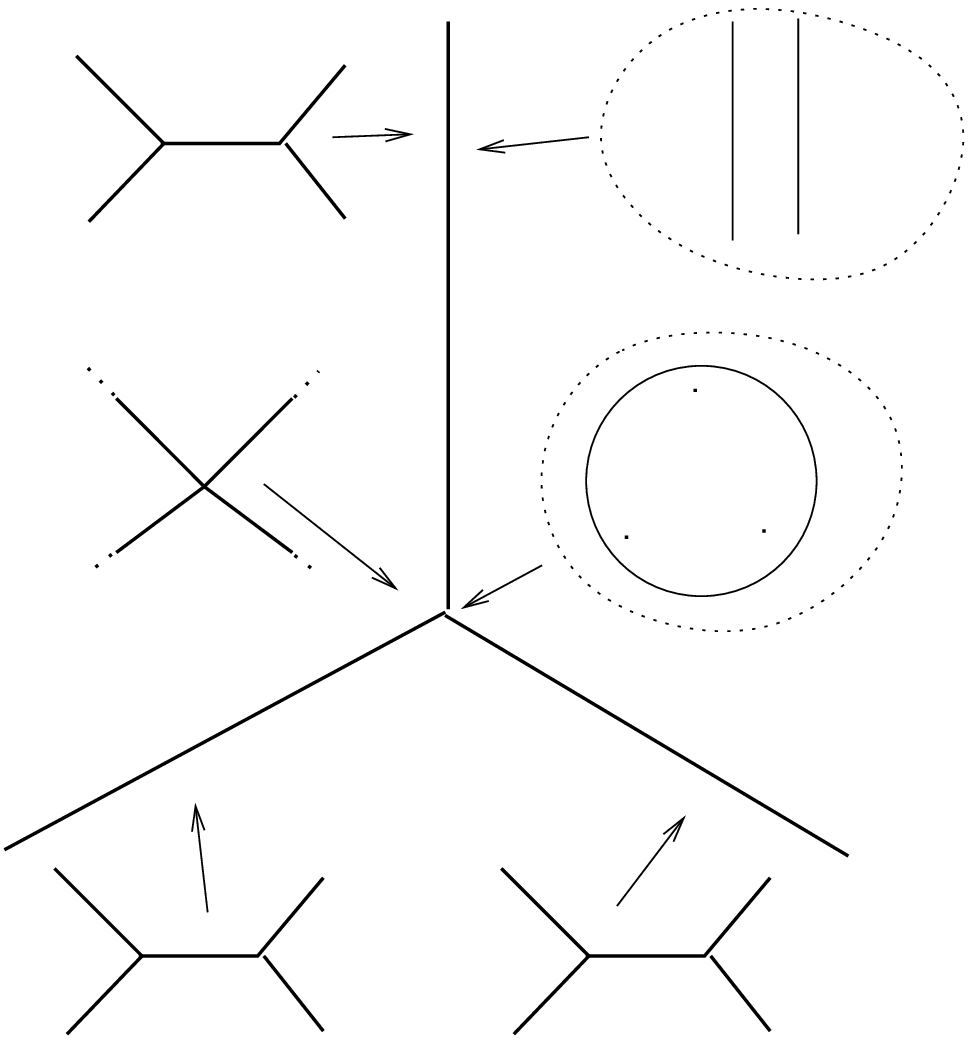}%
\end{picture}%
\setlength{\unitlength}{3947sp}%
\begingroup\makeatletter\ifx\SetFigFont\undefined%
\gdef\SetFigFont#1#2#3#4#5{%
  \reset@font\fontsize{#1}{#2pt}%
  \fontfamily{#3}\fontseries{#4}\fontshape{#5}%
  \selectfont}%
\fi\endgroup%
\begin{picture}(4852,5100)(3729,-6001)
\put(3901,-1036){\makebox(0,0)[lb]{\smash{{\SetFigFont{12}{14.4}{\rmdefault}{\mddefault}{\updefault}{ 1}%
}}}}
\put(3961,-2191){\makebox(0,0)[lb]{\smash{{\SetFigFont{12}{14.4}{\rmdefault}{\mddefault}{\updefault}{ 2}%
}}}}
\put(5536,-1126){\makebox(0,0)[lb]{\smash{{\SetFigFont{12}{14.4}{\rmdefault}{\mddefault}{\updefault}{ 3}%
}}}}
\put(5491,-2086){\makebox(0,0)[lb]{\smash{{\SetFigFont{12}{14.4}{\rmdefault}{\mddefault}{\updefault}{ 4}%
}}}}
\put(3886,-5116){\makebox(0,0)[lb]{\smash{{\SetFigFont{12}{14.4}{\rmdefault}{\mddefault}{\updefault}{ 1}%
}}}}
\put(3916,-6001){\makebox(0,0)[lb]{\smash{{\SetFigFont{12}{14.4}{\rmdefault}{\mddefault}{\updefault}{ 3}%
}}}}
\put(5371,-5071){\makebox(0,0)[lb]{\smash{{\SetFigFont{12}{14.4}{\rmdefault}{\mddefault}{\updefault}{ 2}%
}}}}
\put(5296,-5911){\makebox(0,0)[lb]{\smash{{\SetFigFont{12}{14.4}{\rmdefault}{\mddefault}{\updefault}{ 4}%
}}}}
\put(6016,-4951){\makebox(0,0)[lb]{\smash{{\SetFigFont{12}{14.4}{\rmdefault}{\mddefault}{\updefault}{ 1}%
}}}}
\put(6076,-5971){\makebox(0,0)[lb]{\smash{{\SetFigFont{12}{14.4}{\rmdefault}{\mddefault}{\updefault}{ 4}%
}}}}
\put(7501,-5131){\makebox(0,0)[lb]{\smash{{\SetFigFont{12}{14.4}{\rmdefault}{\mddefault}{\updefault}{ 2}%
}}}}
\put(7501,-5851){\makebox(0,0)[lb]{\smash{{\SetFigFont{12}{14.4}{\rmdefault}{\mddefault}{\updefault}{ 3}%
}}}}
\put(8581,-1561){\makebox(0,0)[lb]{\smash{{\SetFigFont{12}{14.4}{\rmdefault}{\mddefault}{\updefault}{ $\mathbb C^*$}%
}}}}
\put(8221,-3301){\makebox(0,0)[lb]{\smash{{\SetFigFont{12}{14.4}{\rmdefault}{\mddefault}{\updefault}{ $\bar{\mathcal M}_{0,4}$}%
}}}}
\end{picture}%

\end{figure}

The identification between the normal neighborhood bundles of boundary strata of $\bar {\mathcal M}_{0,4} $ and the gluing data $\mathbb C^*$ is explained above.

\subsection{Exploded $\mathbb T$ curves in smooth manifolds}

A connected smooth manifold $M$ can be regarded as an exploded $\mathbb T$ fibration with base consisting of a point, and fiber over that point being $M$. If $M$ is given an almost complex structure, we can then consider holomorphic exploded $\mathbb T$ curves in $M$. A holomorphic exploded $\mathbb T$ curve in this setting has the information of a nodal holomorphic curve along with a positive length and a $\mathbb C$ bilinear pairing of tangent spaces assigned to each node. This extra information corresponds to exploding the moduli space instead of compactifying it. The exploded $\mathbb T$ curve is stable if the underlying nodal holomorphic curve is stable.

\subsection{Symplectic sum}

Given symplectic manifolds $M_1^{2n}$ and $M_2^{2n}$ along with symplectic submanifolds,
\[N^{2(n-1)}\hookrightarrow M_1^{2n}\text{ and }N^{2(n-1)}\hookrightarrow M_2^{2n}\] so that the normal bundle of $N$ in $M_1$ is dual to the normal bundle of $N$ in $M_2$, we can construct the symplectic sum of $M_1$ and $M_2$ over $N$, $M_1\#_NM_2$. The question of how to find holomorphic curves in $M_1\#_NM_2$ in terms of holomorphic curves in $M_1$, $M_2$ and the normal bundle to $N$ has been answered in \cite{IP}. This involves a degeneration of complex structure which can be viewed as giving the following exploded $\mathbb T$ fibration.

First choose an almost complex structure on $M_1$ and $M_2$ so that $N$ is an almost complex submanifold, and the normal bundle to $N$ in $M_1$ is  a complex bundle isomorphic to the dual of the normal bundle in $M_2$. $M_1$ and $M_2$ are log smooth $\mathbb T$ spaces with a single substrata $N$. Note that the normal neighborhood bundle $\mathcal N_NM_1$ is equal to the normal bundle of $N$ minus the zero section.
  
  The exploded fibration has a base equal to the interval $[1,2]$ with fiber over the strata $1$ being $M_1$, fiber over the strata $2$ being $M_2$. The fiber over any point in the strata $[1,2]$ is equal to $\mathcal N_NM_2=\mathbb C^*\rtimes N$. This is identified with $\mathcal N_NM_1$ using the identification of the normal bundle of $N$ in $M_1$ with the dual of the normal bundle of $N$ in $M_2$. Multiplication by a vector $c$ in the tangent space to $[1,2]$ corresponds to multiplication by $e^c$ if the fiber is viewed as the normal bundle of $N$ in $M_2$, and to multiplication by $e^{-c}$ when viewed as the normal bundle of $N$ in $M_1$.

\subsection{Relationship with tropical geometry}

An easy example of an exploded $\mathbb T$ fibration is $\mathfrak T\mathbb R^n$, which has a base given by $\mathbb R^n$ and fibers given by $\left(\mathbb C^*\right)^n$.  We will now examine holomorphic exploded $\mathbb T$ curves
\[f:\ex C\longrightarrow \mathfrak T\mathbb R^n\]

At any vertex $c$ in the base $C$ of $\ex C$ , we have a holomorphic map defined on the interior of the fiber $\Sigma_c$.
\[f_c:\Sigma_c\longrightarrow \left(\mathbb C^*\right)^n\]
The fact that this is a log smooth $\mathbb T$ morphism implies that this must extend to a meromorphic map over punctures. In particular, in local coordinates $z$ around a puncture, 
\[f_c(z):=(z^{\alpha_1}g_1(z),\dotsc,z^{\alpha_n}g_n(z))\]
where $\alpha_i\in\mathbb Z$ and $g_i$ are holomorphic $\mathbb C^*$ valued functions which extend over punctures.
The restriction of $f_c$ to the normal neighborhood is then given by
\[f_c(z):=(z^{\alpha_1}g_1(0),\dotsc,z^{\alpha_n}g_n(0))\]
(The coordinate $z$ in this case is the restriction of the above coordinate to the normal neighborhood bundle.)

This then defines the fiber maps on the edge attached to this puncture, so if $e$ is a point on this edge, then
\[f_e(z):=(z^{\alpha_1}g_1(0),\dotsc,z^{\alpha_n}g_n(0))\]
This has implications for the map $f:C\longrightarrow \mathbb R^n$ of bases. It means that this edge viewed as oriented away from the vertex $c$ must travel in the integral direction
$\vec\alpha:=(\alpha_1,\dotsc,\alpha_n)$. If the edge has length $l$, then the image of the edge is $l\vec\alpha$. 

$df$ restricted to this edge gives a map $\mathbb Z\longrightarrow \mathbb Z^n$ given by $k\mapsto k\vec\alpha$. Call this the momentum of the edge, and $\vec\alpha$ the momentum of the edge exiting the vertex. The sum of all momentums exiting a vertex is zero. This is because number of poles of a meromorphic function is equal to the number of zeros. It can be viewed as a balancing or conservation of momentum condition on the base of a holomorphic curve. 

\

\begin{figure}[htb]
 \begin{picture}(0,0)%
\includegraphics{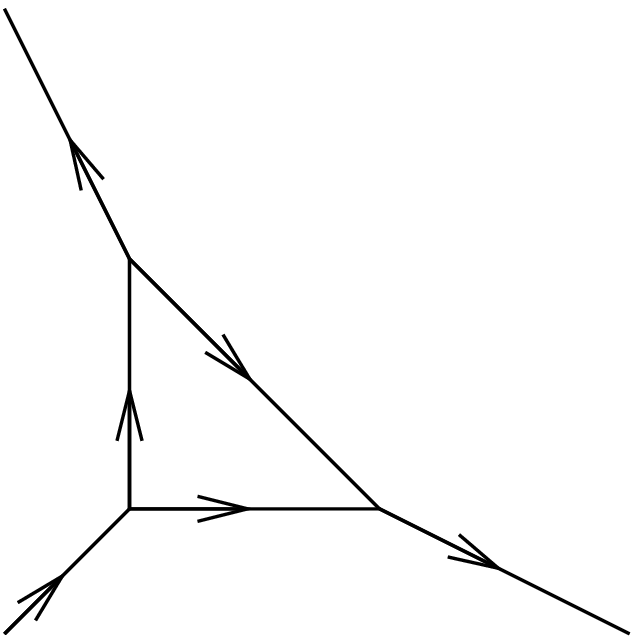}%
\end{picture}%
\setlength{\unitlength}{3947sp}%
\begingroup\makeatletter\ifx\SetFigFont\undefined%
\gdef\SetFigFont#1#2#3#4#5{%
  \reset@font\fontsize{#1}{#2pt}%
  \fontfamily{#3}\fontseries{#4}\fontshape{#5}%
  \selectfont}%
\fi\endgroup%
\begin{picture}(3307,3080)(4179,-3419)
\put(4501,-3361){\makebox(0,0)[lb]{\smash{{\SetFigFont{12}{14.4}{\rmdefault}{\mddefault}{\updefault}{ $(1,1)$}%
}}}}
\put(5101,-3061){\makebox(0,0)[lb]{\smash{{\SetFigFont{12}{14.4}{\rmdefault}{\mddefault}{\updefault}{ $(1,0)$}%
}}}}
\put(5536,-1936){\makebox(0,0)[lb]{\smash{{\SetFigFont{12}{14.4}{\rmdefault}{\mddefault}{\updefault}{ $(1,-1)$}%
}}}}
\put(4501,-541){\makebox(0,0)[lb]{\smash{{\SetFigFont{12}{14.4}{\rmdefault}{\mddefault}{\updefault}{ $(-1,2)$}%
}}}}
\put(6856,-3046){\makebox(0,0)[lb]{\smash{{\SetFigFont{12}{14.4}{\rmdefault}{\mddefault}{\updefault}{ $(2,-1)$}%
}}}}
\put(4186,-2191){\makebox(0,0)[lb]{\smash{{\SetFigFont{12}{14.4}{\rmdefault}{\mddefault}{\updefault}{ $(0,1)$}%
}}}}
\end{picture}%

\end{figure}

The special case when the base of an exploded curve is a trivalent graph and the fiber over all vertices is equal to a three punctured sphere is interesting, because most of the information in the exploded curve can be read off from the base. In this case, the map of bases $f:C\longrightarrow \mathbb R^n$ is called a tropical curve. The technique of counting tropical curves to get holomorphic curve invariants has been used quite successfully (for example, see \cite{Mikhalkin}).

\begin{defn}
A tropical exploded $\mathbb T$ fibration is an exploded $\mathbb T$ fibration in which every fiber is some product $\left(\mathbb C^*\right)^k\times\mathbb CP^{n-k}$ where $\mathbb CP^{n-k}$ is given the log smooth $\mathbb T$ structure with boundary strata given by $(n-k+2)$ generic hyperplanes. 
\end{defn}

The base of a tropical exploded $\mathbb T$ fibration is a tropical space. Holomorphic curves in these tropical exploded $\mathbb T$ fibrations with tropical domains have bases which are tropical curves.

\section{Fiber product}\label{fiber product}

The category of exploded fibrations has a good intersection theory.

\begin{defn}
 Given two exploded morphisms with the same target, 
 \[\ex A\xrightarrow{f}\ex C\xleftarrow{g}\ex B\]
 $f$ and $g$ are transverse if all the following maps are transverse restricted
to the interior of their domains
 \[M_a\xrightarrow{f_a}M_{c}\xleftarrow{f_b} M_b\text{ when }f(a)=c=f(b)\]
 \end{defn}

 \
 
 The fiber product of transverse exploded morphisms is well defined, however it is sometimes not an exploded fibration of the type we have described up until now. It may have strata which are locally modeled on open subsets of \[\{b: \alpha_j\cdot b \geq 0 \}\subset \mathbb R^n\] where $\{\alpha_j\}$ is  some set of  integer vectors. There is a correspondingly more permissive definition of log smooth coordinate charts for fibers.
 
 \begin{thm}
  Suppose that $f$ and $g$ are transverse exploded morphisms, 
   \[\ex A\xrightarrow{f}\ex C\xleftarrow{g}\ex B\]
  so that for any pair of strata $A_j^k$ and $B_i^l$ in the base of $\ex A$ and $\ex B$ so that $M_a\fp{f_a}{g_b}M_b$ is nonempty for some $a\in A$ and $b\in B$,  the subset $\{(a,b): f(a)=g(b)\}\subset A\times B$ is an integral affine space locally modeled on open subsets of $[0,\infty)^n$.
  
  Then there exists an exploded fibration called their fiber product 
  \[\ex A\fp fg\ex B\]
  with maps to $\ex A$ and $\ex B$ so that the following diagram commutes
\begin{displaymath}
\begin{array}{cll}
\ex A\fp fg\ex B & \longrightarrow & \ex A \\ 
 \downarrow & × & \downarrow \\ 
\ex B & \longrightarrow & \ex C
\end{array}
\end{displaymath}
and with the usual universal property that given any commutative diagram
\begin{displaymath}
\begin{array}{lll}
\ex D & \longrightarrow  & \ex A \\ 
\downarrow & × & \downarrow \\ 
\ex B & \longrightarrow & \ex C
\end{array}
\end{displaymath}
there exists a unique morphism $\ex D\longrightarrow \ex A\fp fg\ex B$ so that the following diagram commutes
\begin{displaymath}
\begin{array}{llc}
\ex D & \rightarrow & \ex A \\ 
\downarrow & \searrow & \uparrow \\ 
\ex B & \leftarrow & \ex A\fp fg\ex B
\end{array}
\end{displaymath}

 \end{thm}

 \begin{proof}
 
 First, we check the theorem in the case that $g$ is a map of a point into $\ex C$. In this case, the image of $g$ is some point $c\in C$  and a point $x\in M_c$ in the fiber over $c$. As $x$ is the image of a morphism, it must be in the interior of $M_c$.  We are interested in $f^{-1}(x)$. 
 
First consider $f^{-1}(c)$. The base of $f^{-1}(x)$ will be made up of copies of strata of $f^{-1}(c)$, with each strata of $f^{-1}(c)$ appearing a number of times equal to the number of connected components of $f^{-1}(x)$ over it. Call a strata of $f^{-1}(c)$ empty if $f^{-1}(x)$ is empty over the strata. We will discard such strata.   
 
 Note that for any $a\in A$, $f_a$ is transverse to $x$ on the interior of $M_a$. This property is preserved if we take refinements of $\ex A$. We would like $f_a$ restricted to the interior of any strata of $M_a$ to be transverse to $x$. To achieve this, consider a local refinement of $\ex A$ coming from a subdivision of $A$ close to $a$ so that the non empty strata of $f^{-1}(c)$ are strata in the subdivision. (There may be no global subdivision satisfying this, but it is easy to construct one locally). Note that any morphism $\ex D\longrightarrow\ex A$ contained in $f^{-1}(x)$ lifts locally to a unique morphism to the local refinement of $\ex A$.

   After taking such a refinement, if $a\in f^{-1}(c)$ is in a non empty strata, then $f_a$ is constant in any affine direction. Also $f_a$ restricted to the interior of any strata $S\subset M_a$ is transverse to $x$.   This is because $f_a$ restricted to $\mathcal N_SM$ must be transverse to $x$, but if the interior of $f_a(S)$ intersects $x$, it must be contained in the interior of $M_c$, and in such a case $f_a$ restricted to $\mathcal N_SM_a$ is just the pullback of $f_a$ restricted to $S$. This makes $f_a^{-1}(x)$ some number of log smooth spaces. Affine coordinates are just the restriction of affine coordinates on $M_a$, and boundary defining functions are the restriction of boundary defining functions for strata that $f_a^{-1}(x)$ intersects. A single boundary of $M_a$ may correspond to several boundaries in $f_a^{-1}(x)$. Note that $f_a^{-1}(x)$ does not depend on the particular choice of local refinement.

   If $\iota:a\mapsto \iota(a)$ is an inclusion in $f^{-1}(c)$ and $S$ is the strata in $M_a$ corresponding to $\iota$, then 
  
  \[f_{\iota(a)}^{-1}(x)=\left(f_a\vert_{\mathcal N_SM_a}\right)^{-1}(x)\]

  If we take the non empty strata of $f^{-1}(c)$ to be the base of $f^{-1}(x)$, $f^{-1}_a(x)$ to be the fibers,  and restrict parallel transport and isomorphisms $\iota^\dag$ from (the local refinements of) $\ex A$, then $f^{-1}(x)$ satisfies all the axioms of an exploded fibration apart from having connected fibers. 
  To remedy this, have a copy of each strata of $f^{-1}(c)$ for each fiber component, and connect them by inclusions if one fiber is the normal neighborhood of some strata in the other. This exploded fibration has the required universal property.
  
\
  
  Now, we consider the general case of the transverse intersection of two morphisms, $f$ and $g$.
 
Cover the target $\ex C$ with exploded coordinate charts.
 \[\begin{split}
    &\ex U_\alpha\longrightarrow \ex C
    \\&\downarrow
    \\&\ex R^n
   \end{split}\]

 Note that $f^{-1} \left(\ex U_\alpha\right)$ is an exploded fibration, so we can restrict $f$ to $f^{-1} \left(\ex U_\alpha\right)$.
 
 On these charts, we can consider
 \[f-g:f^{-1}\left(\ex U_\alpha\right)\times g^{-1}\left(\ex U_\alpha\right)\longrightarrow \ex R^n\]
 
 Recall that $\ex R^n$ is the exploded fibration with base equal to $\mathbb R^n$ and fibers all equal to $\mathbb R^n$. As $f$ is transverse to $g$, $(f-g)$ is transverse to the point $0$ in the fiber over $0$, and $(f-g)^{-1}(0)$ is an exploded fibration. This constructs $f^{-1} \left(\ex U_\alpha\right)\fp fg g^{-1} \left(\ex U_\alpha\right)$. This has the correct universal property because any commutative diagram
  \begin{displaymath}
\begin{array}{ccc}
\ex D & \longrightarrow & f^{-1}\left(\ex U_\alpha\right) \\ 
\downarrow &  & \downarrow \\ 
g^{-1}\left(\ex U_\alpha\right) & \longrightarrow & \ex U_\alpha
  \end{array}
  \end{displaymath}
gives a morphism $\ex D\longrightarrow f^{-1}\left(\ex U_\alpha\right)\times g^{-1}\left(\ex U_\alpha\right)$ contained in the inverse image of $0$.

  Doing the same for all coordinate charts and patching together the result define an exploded fibration $\ex A\fp fg\ex B$ with the required universal property.

\end{proof}
 
Note that if $f$ and $g$ are transverse \exploded morphisms satisfying the above condition on the base, $\ex A\fp fg\ex B$ is an \exploded fibration. 

\

An important example is given by the forgetfull morphisms 
\[\pi:\ex M_{g,n+1}\longrightarrow \ex M_{g,n}\]
This map is a surjective submersion restricted to the interior of each fiber, so any smooth exploded morphism $f:\ex F\longrightarrow \ex M_{g,n}$ must be transverse to $\pi$. 

\begin{defn}
 An exploded fibration is complete if any metric on it is complete, and the base has a finite number of strata, all of which are complete with respect to the integral affine connection. (A metric on an exploded fibration $\ex B$ is a symmetric, positive definite, exploded section of $T^*\ex B\otimes T^*\ex B$.)
\end{defn}

\begin{defn}
 An exploded fibration is compact if it is complete and the base is compact.
\end{defn}

\begin{lem}
 If $\ex A$ and $\ex B$ are complete and $f: \ex A\longrightarrow \ex C$ and $g:\ex B\longrightarrow \ex C$ are transverse, then $\ex A\fp fg\ex B$ is complete. It $\ex A$ and $\ex B$ are compact, then $\ex A\fp fg\ex B$ is compact.
\end{lem}

 We now want to say that two exploded morphisms to the same target are generically transverse. In the smooth case, statements of the following sort can be usefull if $f:M\longrightarrow N$ is transverse to $g$ on the boundary of $M$ and $g$ is proper, then there is a small isotopy of $f$ fixed on the boundary to a map transverse to $g$.
 In this exploded setup, instead of a boundary, we will have a subset of strata which are closed under inclusion. This means that if $B_j^k$ is inside this subset, then the target of any inclusion of $B_j^k$ must be contained in this subset.
For isotopys, instead of the unit interval $[0,1]$, we will have the explosion of the unit interval, Expl$[0,1]$.

\begin{lem}\label{transversality}
 Given two smooth exploded morphisms, 
 \[\ex A\xrightarrow{f}\ex C\xleftarrow{g}\ex B\]
  If $\ex A$ and $\ex B$ are complete, and the strata of $A$  are simply connected, then  there exists a small isotopy 
   \[F:\text{Expl}[0,1]\times \ex A\longrightarrow\ex C\]
     so that
      \begin{enumerate}
      \item the restriction to any fiber of $\ex A$, $F_a:[0,1]\times M_a\longrightarrow M_{F(a)}$ is a $C^\infty$ small isotopy.
 \item $F(t,\cdot)$ is independent of $t$ in a neighborhood of $t=0$ and $t=1$,  equal to $f$ when $t$ is small, and transverse to $g$ when $t$ is close to $1$.
      \end{enumerate}

Moreover, if $f$ is transverse to $g$ on a subset of strata closed under
inclusions, then the isotopy can be chosen to be fixed on those strata.    
\end{lem}

\begin{proof} 

We construct the isotopy first over fibers that have no boundary strata and then extend it over fibers with strata of higher codimension. 

If $M_a$ has no boundary strata, it is some affine $\mathbb R^n$ bundle over a compact manifold. As the isotopy can't change $df_a$ in affine directions, the isotopy is determined by its restriction to some section of this bundle. Note that to achieve transversality, we can quotient the target of $f_a$ by the directions spanned by $df$ of affine directions. The task is then to find an isotopy so that the image of this section is transverse to the image of $g$. As $\ex B$ is complete, the image of $g$ in this fiber consists of the image of a finite number of smooth proper maps of manifolds with boundary and corners, so a generic smooth map is transverse to $g$ and there exists a $C^\infty$ small isotopy of $f_a$ making it transverse to $g$ as required. (We chose the trivial isotopy if $f_a$ is already transverse.)

This isotopy can be extended compatibly via parallel transport to all fibers of $\ex A$ that have no boundary (this uses that the strata of $A$ are simply connected, so there is no monodromy).
Now suppose that the isotopy has been defined compatibly on all fibers that only have strata of codimension less than $k$. We now must extend it to an isotopy on fibers with strata of maximum codimension $k$. 

Suppose that $M_a$ is such a fiber. As any strata $S$ of nonzero codimension has substrata of dimension at most $k-1$, we have defined the isotopy already on $\mathcal N_SM_a$. Viewing our isotopy as given by exponentiating some smooth, time dependant vector field, we just need to extend this vectorfield over $M_a$. It is defined on the boundary of $M_a$ so that it is continuous, and smooth restricted to any boundary or corner strata. As such it extends to a smooth vectorfield on the interior of $M_a$, which can be chosen $C^\infty$ small depending on its $C^\infty$  size on the boundary strata. (One can see locally why this is true: A continuous function defined the boundary of $[0,\infty)^n$ which is smooth on boundary and corner strata extends smoothly to the interior as follows. Denote by $\pi_S$ the projection to the strata $S$, and let $\sigma(S)$ be 1 if $S$ has odd codimension, and $-1$ if $S$ has even codimension. Then we can set $f:=\sum\sigma(S)f\circ\pi_S$.) 

 We now have a $C^\infty$ small isotopy of $f_a$ that restricts to any nontrivial normal neighborhood bundle to be our previously defined isotopy. As this restriction to the normal neighborhood bundle is transverse to the image of $g$ near $t=1$, and the space of maps transverse to the image of $g$ is open, this isotopy is transverse to the image of $g$ near the boundary and $t=1$. We can then further modify this isotopy on the interior of $M_a$ to be transverse to $g$ everywhere near $t=1$.  This isotopy can be chosen $C^\infty$ small, and if $f_a$ is already transverse to $g$ and we have chosen the isotopy trivial on the normal neighborhood bundles of  boundary strata of $M_a$, the this isotopy can be trivial. 

Thus we can extend the isotopy as required.

\end{proof} 

A very similar statement holds for exploded $\mathbb T$ fibrations, however there is one small complication: Expl$[0,1]\times\ex A$ is not an \exploded fibration. We can define a category of mixed exploded fibrations which includes fiber products of exploded torus fibrations and exploded fibrations.
 The fibers of these are locally modeled on $\{z^{\alpha_i}\in [0,\infty)\}\subset \mathbb C^n\times\mathbb  R^k$ where $\{\alpha_i\}$ are multiindices in $\mathbb Z^n$ which form a basis for some subspace of $\mathbb Z^n$, and log smooth functions are the restrictions of functions of the form $z^\alpha g$ where $g$ is a smooth $\mathbb C^*$ valued function. 
 The case when the set $\{\alpha_i\}$ spans $\mathbb Z^n$ corresponds to exploded fibrations, and the case where this set is empty corresponds to exploded torus fibrations. Definitions of exploded fibrations using this mixed version of log smooth spaces are the same as our earlier definitions. The results of section \ref{refinements} and this section apply for this more general category.

\section{Moduli space of smooth exploded $\mathbb T$ curves}\label{perturbation theory}

In order to describe the structure on the moduli space of smooth exploded curves, we will need to use notions from the theory of polyfolds developed by Hofer, Wysocki and Zehnder in \cite{polyfold}.

\begin{defn}
 A Frechet splicing core $K$ is the fixed point set of a smooth family of projections
 \[\pi:F\times[0,\infty)^k\rightarrow F\times[0,\infty)^k\]
 where $\pi(\cdot,x)$ gives a projection on the Frechet space $F\longrightarrow F$ and $\pi(\cdot,0)$ is the identity. 

\end{defn}

A log smooth function  $f\in\logC^\infty\left(F\times [0,\infty)^k\right)$ is a function defined on the interior of $F\times[0,\infty)^k$ so that 
\[f=\sum \alpha_i\log x_i+g\]
where $\alpha_i \in\mathbb Z$ and $g$ is a smooth function on $F\times[0,\infty)^k$.

\begin{defn}
 A log smooth Frechet polyfold is a Hausdorff second countable topological space $M$ along with a sheaf of log smooth functions $\logC^\infty(M)$ so that around each point $p\in M$, there exists an open set $U_p$ and a homeomorphism onto some Frechet splicing core $\phi_p :U_p\longrightarrow K\subset F_p\times[0,1)^{k_p}$ sending $p$ to the fiber over $0$ so that 
 \[\logC^\infty(U_p)=\phi_p^*\left(\logC^\infty(F_p\times[0,\infty)^{k_p})\right)\]
$p$ is said to be in a strata of codimension $k_p$. 
\end{defn}

\

Note that if $M$ is finite dimensional it is simply a log smooth space, as then splicing cores must be diffeomorphic to $F\times[0,\infty)^k$.

The subset of $M$ consisting of points with codimension  $k$ is a smooth Frechet manifold. For each connected component of these points of codimension $k$, there exists a unique log smooth Frechet polyfold $S$ and a continuous map $S\longrightarrow M$ so that the interior of $S$ maps diffeomorphically to this component. Call all such $S$ the strata of $M$.
 
We can define affine $\mathbb R^n$ bundles over log smooth Frechet polyfolds, log smooth maps, normal neighborhood bundles and log smooth morphisms exactly as in the finite dimensional case. With this done, the definition of a Frechet exploded fibration is identical to the definition of an exploded fibration where fibers are  log smooth Frechet polyfolds.

We can similarly define Frechet exploded $\mathbb T$ fibrations.

\begin{defn}
 Given a closed $2$ form $\omega$ on the exploded $\mathbb T$ fibration $\ex B$, the smooth exploded $\mathbb T$ curve
 $f:\ex C\longrightarrow\ex B$ 
 is $\omega$-stable if $C$ has a finite number of strata and the integral of $f^*\omega$ is positive on any unstable fiber $\Sigma_v$ over a vertex of $C$. 
\end{defn}

\

What follows is a rough sketch of the perturbation theory of $\omega$-stable 
holomorphic curves.

\begin{conj}
 The moduli space $\ex M^{sm}_{g,n}(\ex B,\omega)$ of $\omega$-stable smooth exploded curves in $\ex B$ with genus $g$ and $n$ marked points has a natural structure of a Frechet orbifold exploded fibration.
\end{conj}

If $\ex B$ has an almost complex structure $J$, we can take $\dbar$ of a smooth exploded $\mathbb T$ curve $f:\ex C\longrightarrow \ex B$. If $j$ denotes the complex structure on $\ex C$, this is defined as 
\[\dbar f:=\frac 12 (df+J\circ df\circ j)\]
This takes values in antiholomorphic sections of $T^*\ex C\otimes f^*(T\ex B)$ which vanish on the fibers over the edges of $\ex C$.

\begin{conj}
The space of these antiholomorphic sections has the structure of a Frechet orbifold exploded fibration $\ex E_{g,n}(\ex B,\omega)$, which is a vector bundle over $\ex M^{sm}_{g,n}(\ex B,\omega)$.
The $\dbar$ equation defines a smooth section of this vector bundle.
\[\begin{split}&\ex E_{g,n}(\ex B,\omega)
   \\& \downarrow \uparrow\dbar
   \\&\ex M_{g,n}^{sm}(\ex B,\omega)
  \end{split}
\]
\end{conj}

\begin{conj}\label{Fredholm}
The section $\dbar$ is `Fredholm' in the following sense:  There is a cover of
$ \ex M_{g,n}(\ex B,\omega)$ by exploded coordinate neighborhoods $(\ex U, G)$
covered by coordinates $(\ex E,G )$ on $\ex E_{g,n}(\ex B,\omega)$ so that there exists the following
\begin{enumerate}
\item An identification of $\ex U$ with a coordinate neighborhood in a vector bundle over a finite dimensional exploded $\mathbb T$ fibration $\ex F$
\[\begin{split}
\ex U\xrightarrow{\psi} &\ex V
\\&\downarrow
\\&\ex F   
  \end{split}
\]
\item An identification of $\ex E$ with a coordinate neighborhood in a fibration 
\[\begin{split}\ex E\subset \ex V &\oplus\ex{\tilde E}
   \\&\downarrow
   \\&\ex F
  \end{split}
\]
so that the projection $\ex E\rightarrow \ex U$ is given by $(v,e)\mapsto v $, and the zero section $\ex U\longrightarrow \ex E$ is given by $v\mapsto (v,0)$.

\item An exact sequence of vectorbundles over $\ex F$
\[0\longrightarrow \mathbb R^k\times\ex F\longrightarrow \ex{\tilde E}\xrightarrow{\pi_{\ex V}} \ex V\longrightarrow 0\]
so that the graph of $\pi_{\ex V}\circ\dbar$ is the diagonal in $\ex V\oplus \ex V$
\end{enumerate}
\end{conj}

In the non algebraic case, the regularity of maps and identifications in the above conjecture may be questionable. 

\

The above conjecture says roughly that the intersection theory of $\dbar$ with the zero section can be modelled locally on the following finite dimensional problem: first, the solutions of the equation $\pi_{\ex V}\circ\dbar$ can be identified  with $\ex F$. $\dbar$ of such a solution then has values in the $\mathbb R^k$ from the above exact sequence, so we get a map $\dbar: \ex F\longrightarrow \mathbb R^k$. We are interested in the transverse intersection of this map with $0$.

\

We would like to take the transverse intersection of $\dbar$ with the zero section and obtain the moduli space of holomorphic curves as a smooth exploded $\mathbb T$ fibration. The problem with this is that as we are dealing with orbifolds, we can't necessarily perturb the zero section to get it transverse to the zero section. Following \cite{Orbifolds}, we can take a resolution of the zero section by a nonsingular weighted branched exploded fibration $\ex{\tilde M}\rightarrow \ex M^{sm}_{g,n}(\ex B,\omega)$, and then perturb that to be transverse to $\dbar$, obtaining a representation of the virtual moduli space of holomorphic curves inside $\ex M^{sm}_{g,n}(\ex B,\omega)$ as a smooth weighted branched exploded $\mathbb T$ fibration.

\

{\bf Acknowledgments:} I wish to thank Yakov Eliashberg, Pierre Albin, Tyler Lawson, Kobi Kremnizer, Tomasz Mrowka, and Eleny Ionel for helpfull conversations and suggestions which aided the development of exploded fibrations.

\end{document}